\documentclass[reqno,11pt]{amsart}
\usepackage{xcolor}
\usepackage[margin=1in]{geometry}
\usepackage{amsthm,amsmath,amssymb,dsfont}
\usepackage{mathrsfs,amsfonts,functan,extarrows,mathtools}

\usepackage[colorlinks,
linkcolor=red,
citecolor=blue
]{hyperref}
\usepackage{indentfirst, latexsym, enumerate,graphicx}
\usepackage{float}
\usepackage{relsize}
\usepackage{marginnote}
\usepackage{stmaryrd}
\usepackage{esint}
\usepackage{bm}
\usepackage{caption}
\usepackage{cite}
\usepackage{color}
\usepackage{subfigure}
\usepackage{paralist}
\allowdisplaybreaks 

\theoremstyle{plain}
\newtheorem{theorem}{Theorem}[section]

\newtheorem{lemma}{Lemma}[section]
\newtheorem{proposition}{Proposition}[section]
\newtheorem{remark}{Remark}[section]

\numberwithin{equation}{section}

\usepackage{appendix}
\usepackage{xcolor}

\def\rmd{\mathrm{d}}
\def\P{\mathbf{P}}
\def\Pp{\mathbf{P}^\perp}

\begin{document}

\title[Optimal time decay of the Boltzmann equation with cutoff soft potentials]
{Optimal time decay of the Boltzmann equation with cutoff soft potentials in the whole space}

\author[Y.J.  Lei]{Yuanjie Lei}  
\address[YJL]{School of Mathematics and Statistics,  Huazhong University of Science and Technology,  and the Hubei Key Laboratory of Engineering Modeling and Scientific Computing,  Wuhan  430074, P. R. China}
\email{leiyuanjie@hust.edu.cn}
 
\author[X.J.  Wei]{Xiongju Wei}   
\address[XJW]{School of Mathematics and Statistics,  Huazhong University of Science and Technology,  and the Hubei Key Laboratory of Engineering Modeling and Scientific Computing,  Wuhan  430074, P. R. China}
\email{weixiongju@hust.edu.cn}
    
\author[T. Yang]{Tong Yang}  
\address[TY]{Department of Applied Mathematics, The Hong Kong Polytechnic University,
Kowloon, Hong Kong, P. R. China}
\email{t.yang@polyu.edu.hk}
    
\author[H.J. Zhao]{Huijiang Zhao}
\address[HJZ]{School of Mathematics and Statistics, Wuhan University, Wuhan 430072, P. R. China; Computational Science Hubei Key Laboratory, Wuhan University, Wuhan 430072, P. R. China}
\email{hhjjzhao@whu.edu.cn}   

\begin{abstract}
    Since Strain's work \cite{Strain-2012-KRM} on almost optimal time decay rates for the microscopic part of perturbative classical solutions to the Boltzmann system with soft potentials via the time-velocity splitting method, how to obtain the optimal decay of the microscopic part has remained unsolved. In this paper, we resolve this problem by developing a nonlinear energy estimate method different from the time-velocity splitting approach.
    The proof relies on a combination of a hierarchy of time-weighted energy estimates, the damping factor \(e^{-\nu(v)(t-\tau)}\), and a dyadic decomposition of velocity space.
\end{abstract}

\keywords{Boltzmann equation, optimal time decay, soft potentials}

\subjclass[2020]{35A01,  35B35, 76P05}

\maketitle

\setcounter{equation}{0}
\indent \allowdisplaybreaks

\tableofcontents

\section{Introduction}
    \subsection{Problem}
    We consider the Boltzmann equation
    \begin{equation}\label{BQ-1}
    	\partial_t F+v\cdot\nabla_xF=Q(F,F),
    \end{equation}
    with the initial data $F(0,x,v)=F_0(x,v)$. 
    Here $F(t,x,v)\geq0$ is the distribution function for the particles at time $t\geq 0$, with position $x = (x_1, x_2, x_3)\in \mathbb{R}^3$ and the velocity $v = (v_1, v_2, v_3)\in \mathbb{R}^3$.
    The collision operator $Q(F,F)$ is defined by
    \begin{equation*}
    	Q(F,F)
    	=\int_{\mathbb{R}^3\times\mathbb{S}^2} B(|v-u|,\omega)
        [ F(u')F(v')-F(u)F(v) ] \rmd u\rmd \omega,
    \end{equation*}
    where $v, u$ and $v^{\prime},u^{\prime}$ are the velocities of a pair of particles before and after collision. They are connected through the formulas
	\begin{align*} 
		v^{\prime}=v-[(v-u)\cdot\omega]\omega,\quad u^{\prime}=u+[(v-u)\cdot\omega]\omega, \quad \omega \in \mathbb{S}^2,
	\end{align*}
    and satisfy the conservation laws of momentum and energy
    \begin{align*}
        v+u=v'+u',
        \quad
        |v|^2+|u|^2=|v'|^2+|u'|^2.
    \end{align*}
    The collision kernel $B(|v-u|,\omega)$ satisfies
    \begin{align*}
        B(|v-u|,\omega)=|v-u|^\gamma b(\cos\theta), \quad \cos \theta=\frac{|(v-u)\cdot\omega|}{|v-u|}, \quad \theta\in\left[0,\frac\pi2\right].
    \end{align*}
    We consider the soft potential case $-3<\gamma<0$ in this paper. 
    Throughout the paper, we assume that $b(\cos\theta)$ satisfies the Grad's angular cutoff assumption (cf. \cite{Grad-1958,Cercignani-Illner-Pulvirenti-1994})
    \begin{equation}
    	0 \leq b(\cos\theta) \leq C |\cos \theta |. \notag
    \end{equation}
    
    For the Boltzmann equation with Grad's angular cutoff cases near Maxwellians, in 1974,
    Ukai \cite{Ukai-1974-PJA} first gave the spectral analysis of Boltzmann equation for hard potentials in a periodic box and constructed the global-in-time solutions; subsequently, Ukai-Yang \cite{Ukai-Yang-2006-AA} extended the analysis to a broader class of solution spaces; see also \cite{Ukai-1976, Zhong-2014-SCM}.  
    Guo \cite{Guo-2003-ARMA} constructed the global classical solutions for all soft potentials with $-3<\gamma<0$ in a periodic box. 
    In the whole space, Guo \cite{Guo-2004-IUMJ} and Liu-Yang-Yu \cite{Liu-Yang-Yu-2004-PD} independently developed approaches based on the macro–micro decomposition and nonlinear energy estimates, providing a framework for studying global solutions and their large-time behavior.

   	A central question in this theory is how rapidly perturbative solutions converge to the global Maxwellian as time tends to infinity. For the hard-potential case, Ukai's spectral approach \cite{Ukai-1974-PJA} and the compensating-function method developed by Kawashima \cite{Kawashima-1990-JJAM} provide important tools for obtaining optimal decay rates. These ideas have also been developed for more general Boltzmann-type systems; see \cite{Yang-Yu-2011-CMP, Duan-Strain-2011-ARMA, Duan-Strain-2011-CPAM} and the references therein.
    Based on the optimal time decay rates such as 
    \[
        \|\nabla^k_x f\|^2\sim (1+t)^{-\frac32-k},\quad k=0,1,
    \]
    since the microscopic part satisfies 
    \begin{equation}
        \partial_t {\bf P}^{\perp}f+v\cdot \nabla_x\Pp f +L  {\bf P}^{\perp}f
        = - v\cdot \nabla_x{\bf P} f+ \P(v\cdot\nabla_x f)+\Gamma(f,f), \notag
    \end{equation}
    one has for some $\lambda >0$,
    \begin{equation}
        \|{\bf P}^{\perp}f \|^2
        \lesssim 
        e^{-\lambda t} \|{\bf P}^{\perp}f _0\|^2 
        +\int_0^t e^{-\lambda (t-\tau)}\|{\nabla_x }f \|^2\rmd \tau+\dots
        \lesssim (1+t)^{-\frac{5}{2}}, \nonumber
    \end{equation}
    which is the optimal time decay rate of the microscopic part in the whole space.
    
    However, since the linearized operator $L$ has no spectral gap for soft potentials, this gives rise to more analytic difficulties than the study of hard potentials. 
    Recently, Yang-Yu\cite{Yang-Yu-2024-IUMJ} analyzed the spectrum structure of the Boltzmann equation for the whole range of soft potentials, where they obtained the optimal time decay rate as follows,
    \begin{align}
        \| f \|_{N,\ell,q,\theta}
        =\left\| \langle v \rangle^\ell e^{\frac{q}{4}\langle v \rangle^\theta} \|f(t,x,v)\|_{H^N_x(\mathbb{R}^3)} \right\|_{L^\infty_v(\mathbb{R}^3)}
        \sim (1+t)^{-\frac34}. \notag
    \end{align}
    In addition, Guo–Wang \cite{Guo-Wang-2012-CPDE} developed a general energy method via negative Sobolev space interpolation for optimal decay rates of dissipative equations in the whole space, while Sohinger–Strain \cite{Sohinger-Strain-2014-Adv} introduced a homogeneous Besov–Lipschitz interpolation method for optimal decay estimates of the Boltzmann equation with hard and soft potentials.
    Compared to the macroscopic part, to obtain the faster time decay rates of the microscopic part of the solution, by using the compensation function, interpolation and time-velocity splitting method, Strain \cite{Strain-2012-KRM} showed that for the soft potentials
    \begin{equation}\label{decay-micro}
        \|{\bf P}^{\perp}f \|^2
        \lesssim (1+t)^{-\frac52+\varepsilon},
    \end{equation}
    where $\varepsilon >0$ is a constant. 
    A natural question is whether the time decay estimate \eqref{decay-micro} holds when $\varepsilon=0$. 
    Therefore, we aim to study the optimal temporal decay estimate for the microscopic part of solutions of the Boltzmann equation, which is the main purpose of this paper.

    \subsection{Reformulation}
    We begin by reformulating the problem. 
    The normalized global Maxwellian $\mu$ is defined as
    \[
        \mu(v)=(2\pi)^{-\frac32}e^{-\frac{|v|^2}2}.
    \]
    Set the perturbation $f=f(t,x,v)$ as
    \[
        F=\mu+\mu^{\frac12}f.
    \]
    Then the Boltzmann equation \eqref{BQ-1} can be rewritten as
    \begin{equation}\label{f-perturbation}
    	\partial_t f +v\cdot \nabla_x f+Lf=\Gamma(f,f),
    \end{equation}
    with the initial data
    \begin{align}\label{f0}
        f(0,x,v)=f_0(x,v),
    \end{align}
    where the linearized collision operator $L$ is given by
    \begin{equation*}
    	Lf = -\mu^{-\frac12} \left\{{Q (\mu,\mu^{\frac12}f )+ Q (\mu^{\frac12}f,\mu )} \right\},
    \end{equation*}
    and the nonlinear collision term ${\Gamma}(f,f)$ is given by
    \begin{equation}\label{Gamma-def}
    	{\Gamma}(f,f) =\mu^{-\frac12}Q (\mu^{\frac12}f,\mu^{\frac 12}f ).
    \end{equation}
    
    The linearized operator $L$ can be decomposed as
    \begin{align*}
        L=\nu-K,
    \end{align*}
    where the collision frequency $\nu(v)$ is defined by
    \begin{align}
        \nu(v)=\int_{\mathbb{R}^3\times\mathbb{S}^2}|v - u|^\gamma b(\cos\theta)\mu(u){\rm d}u{\rm d}\omega,\notag
    \end{align}
    and satisfies
    \begin{align*}
        \nu(v) \sim (1 + |v|)^{\gamma},\quad -3<\gamma<0.
    \end{align*}
    The operator $K=K_2-K_1$ is a compact operator on $L^2(\mathbb R^3_{v})$ and defined by
    \begin{align}
        K_1 f(v)=&\int_{\mathbb R^3\times\mathbb S^2} |v-u|^\gamma b(\cos\theta) \mu^\frac12(u)\mu^\frac12(v)f(u){\rm d}u{\rm d}\omega,  \nonumber\\
        K_2 f(v)=&\int_{\mathbb R^3\times\mathbb S^2} |v-u|^\gamma b(\cos\theta)\mu^\frac12(u)\left\{\mu^\frac12(u^\prime)f(v^\prime)+\mu^\frac12(v^\prime)f(u^\prime)\right\}{\rm d}u{\rm d}\omega . \nonumber
    \end{align}
    
    As in \cite{Guo-2004-IUMJ, Glassey-1996}, the linearized collision operator $L$ is non-negative. The null space of ${L}$ is the five-dimensional space defined by
    \begin{equation*}
    	{\mathcal{ N}}={\textrm{Span}}\left\{\mu^{\frac12}, ~ v_i\mu^{\frac12} ,(| v|^2-3 )\mu^{\frac12}\right\},\quad 1\leq i\leq3.
    \end{equation*}
    For any function $f(t,x,v)$, we define ${\bf P}$ as the orthogonal projection operator from $L^2({\mathbb{R}}^3_v)$ to the null space $\mathcal{N}$. Consequently, for any given function $f(t,x,v)\in L^2(\mathbb{R}^3_v)$, we can express it as
    \begin{align}
        {\bf P}f =\{ a+b\cdot v+c(| v|^2-3)\}\mu^{\frac12},\notag
    \end{align}
    where
    \begin{equation}\nonumber
        \left\{
        \begin{aligned}
            a(t, x)&= \int_{\mathbb{R}^3_v}f \mu^{\frac12}\rmd v,\\
            b_i(t, x)&= \int_{\mathbb{R}^3_v} f v_i \mu^{\frac 12}\rmd v,\quad i =1,2,3,\\
            c(t, x)&= \frac{1}{6}\int_{{\mathbb{R}}^3_v} f (|v|^2-3) \mu^{\frac12} \rmd v.
        \end{aligned}
        \right.
    \end{equation}
    Hence we obtain the unique macro-micro decomposition
    \begin{equation}
    	f(t,x, v)={\bf P}f(t,x, v)+{\bf P^{\perp}}f(t, x, v),\notag
    \end{equation}
    where ${\bf P}f$ and ${\bf P}^\perp f=\{ {\bf I}- \P \}f$ are referred to as the macroscopic and microscopic components of $f(t,x,v)$, respectively. For a corresponding macro-microscopic decomposition applied to the Boltzmann equation, see \cite{Guo-2004-IUMJ, Liu-Yang-Yu-2004-PD}. 

    We introduce the following dyadic decomposition of the velocity space \(\mathbb{R}^3_v\):
    \begin{align}
        \mathbb{R}^3_v=\mathbb{V}_0 \cup \Big(\bigcup_{j\ge1}\mathbb V_j\Big),\label{dyadic}
    \end{align}
    with
    \[
        \mathbb{V}_0=\{v\in\mathbb{R}^3 \mid |v|< 2\},
    \]
    and
    \[
        \mathbb{V}_{j}=\{v\in\mathbb{R}^3 \mid 2^{j}\leq |v| < 2^{j+1},\ j\in\mathbb{Z}^+\}.
    \]
    
   Furthermore, we introduce more notation. For $\ell\geq0$,  define the velocity weight function $w_\ell$ by
    \begin{align}
        w_{\ell}=w_{\ell}(v)=\langle v \rangle^{-\gamma \ell}, \notag
    \end{align}
    with $\langle v \rangle=\sqrt{1+|v|^2}$. 
    Define the instantaneous energy functional $\mathcal{E}_{N,\ell}(t)$
    for $\ell \geq N$ by
    \begin{align}
        \mathcal{E}_{N,\ell}(t) &\sim \sum_{|\alpha|+|\beta|\leq N} \|w_{\ell-|\beta|}\partial^\alpha_x\partial^\beta_vf(t)\|^2, \label{energy functional}
    \end{align}
    and define the dissipation rate functional $\mathcal{D}_{N,\ell}(t)$ by
    \begin{align}\label{dissipation functional}
        \mathcal{D}_{N,\ell} (t)
        \sim \sum_{1\leq |\alpha| \leq N} \|\partial_x^\alpha \P f(t)\|^2
        +\sum_{|\alpha|+|\beta| \leq N} \|w_{\ell-|\beta|} \partial_x^\alpha \partial_v^\beta \Pp f(t) \|^2_\nu,
    \end{align}
    where $N$ is an integer.

    \subsection{Main result}
    We now state our main optimal time decay results for the Boltzmann equation \eqref{f-perturbation} and \eqref{f0}.
    \begin{theorem}\label{thm}
        Let $-3<\gamma<0$, $N\geq4$, $\ell_0\geq N$, $\ell\geq \ell_0+\frac 3 2$. Let $f(t,x,v)$ be the solution to the Cauchy problem \eqref{f-perturbation} and \eqref{f0}. Suppose that
		\begin{align}
			\epsilon_0^2\equiv \sum_{|\alpha|+|\beta|\leq N} \|w_{\ell-|\beta|}\partial^\alpha_x\partial^\beta_vf_0(x,v)\|_{L^2_{x,v}}^2+ \| w_{\frac N2} f_0(x,v) \|^2_{L^2_vL^1_x} \notag
		\end{align}
        is sufficiently small. Then for any $t\geq0$ we have the following time decay estimate
        \begin{align}\label{optimal-decay}
            \|w_{\ell_0}\mathbf{P}^{\perp}f(t,x,v)\|_{L^2_{x,v}}^2  \lesssim (1+t)^{-\frac52} \epsilon_0^2.
        \end{align}
    \end{theorem}
    
    \begin{remark}
       In contrast to the result in \cite{Strain-2012-KRM}, which gives
        \[
            \|\mathbf P^\perp f(t)\|^2
            \lesssim (1+t)^{-\frac52+\varepsilon},
            \quad \varepsilon>0,
        \]
        we obtain the optimal decay rate
        \[
            \|w_{\ell_0}\mathbf P^\perp f(t)\|^2
            \lesssim (1+t)^{-\frac52},
        \]
        without any loss in the time decay exponent. The main reason for this improvement is that we avoid the time-velocity splitting used in \cite{Strain-2012-KRM}. Instead, the damping factor $e^{-\nu(v)(t-\tau)}$ in the Duhamel representation is exploited directly, in combination with a time-weighted Lyapunov hierarchy and a dyadic decomposition of velocity space.
    \end{remark}
    \begin{remark}
        Our research strategy may also be applicable to other kinetic models, such as the Vlasov-Poisson-Boltzmann system, the Vlasov-Maxwell-Boltzmann system, and related systems.
    \end{remark}
    \subsection{Strategy of the proof}
    We first recall Strain's approach \cite{Strain-2012-KRM}, which provides the starting point for our analysis. 
    \subsubsection{Strain's result for the soft potentials}
    Strain \cite{Strain-2012-KRM} established weighted time-frequency Lyapunov inequalities for the linearized equation. The Fourier variable $\xi$ is split into the low-frequency regime $|\xi|\le 1$ and the high-frequency regime $|\xi|>1$: the analysis combines microscopic coercivity with Kawashima-type macroscopic interaction functionals in the low-frequency regime, while using the weighted $L^2$ dissipation directly in the high-frequency regime.
    This yields a functional 
    \begin{align}
        \mathcal{E}_\ell(t,\xi) \sim |w_\ell \hat f(t,\xi)|^2_{L^2_v}, \notag
    \end{align}
    satisfying
    \begin{align}
        \partial_t \mathcal{E}_\ell(t,\xi) + \lambda \min\{1,|\xi|^2\} |w_\ell \hat f(t,\xi)|^2_{L^2_\nu} \leq 0.\notag
    \end{align}
    By using interpolation with a family of higher weight functions, the linear decay was obtained. 
    
    \textbf{Time-velocity splitting in \cite{Strain-Guo-2008-ARMA, Strain-2012-KRM}:} 
    for the nonlinear problem, they introduced the following time-velocity splitting:
    \begin{align}
        E(t) = \{\langle v \rangle^{-\gamma-2s} \leq t^{1-p}\}, \quad E^c(t) = \{\langle v \rangle^{-\gamma-2s} > t^{1-p}\},\notag
    \end{align}
    where the low-weight region $E(t)$ is controlled by dissipation, while the high-weight region $E^c(t)$ is handled by extra polynomial velocity weights and smallness of the initial data. Combining the linear decay with the nonlinear energy Lyapunov inequalities and a bootstrap argument gives
    \begin{align}
        \mathcal E_{N,\ell}(t)\lesssim (1+t)^{-3/2}, \quad  \mathcal E^h_{N,\ell}(t)\lesssim (1+t)^{-5/2+\varepsilon},\notag
    \end{align}
    where $\mathcal E^h_{N,\ell}(t)$ contains the microscopic part $\|\Pp f(t)\|^2$.
    Indeed, according to the time-velocity splitting argument as in \cite{Strain-2012-KRM}, one has
    \begin{align}
        |\Pp \hat f(t,\xi)|^2_{L^2_v} \lesssim e^{-\lambda t^p} |\Pp \hat f_0(\xi)|^2_{L^2_v}+\int_0^t e^{-\lambda(t^p-\tau^p)} \left( |\xi|^2|\P \hat f|^2_{L^2_v} + \cdots \right)\rmd \tau, \notag
    \end{align}
    which implies
    \begin{align}
        \| \Pp f(t)\|_{L^2_{x,v}}^2 \lesssim (1+t)^{-\frac52+\varepsilon(p)},\notag
    \end{align}
    where $\varepsilon(p)=1-p>0$ is a small constant.

    \subsubsection{A method different from the time-velocity splitting method }Instead of using the time-velocity splitting method of \cite{Strain-2012-KRM}, by  decomposing the linear operator as $L=\nu-K$ and exploiting the collision frequency $\nu(v)$, we obtain
    \begin{align}
        &\frac{1}{2}\partial_t |w_{{\ell_0} }\Pp f|^2
        +\frac{1}{2}\nu(v) |w_{{\ell_0}}\Pp f|^2
        +\frac{1}{2}\nu(v) |w_{{\ell_0} }\Pp f|^2 \notag\\
        =&- (v\cdot\nabla_x\Pp f) w_{{\ell_0}}^2 \Pp f
        - (v\cdot\nabla_x\P f) w_{{\ell_0} }^2 \Pp f
        + \P(v\cdot\nabla_x f) w_{{\ell_0} }^2 \Pp f \notag\\
        &+ K\Pp f  w_{{\ell_0} }^2 \Pp f
        + \Gamma(f,f) w_{{\ell_0} }^2 \Pp f,\notag
    \end{align}
    which yields
    \begin{align}
        & \|w_{{\ell_0} }\Pp f\|^2
        +\int_0^t \int_{\mathbb{R}_v^3\times\mathbb{R}_x^3} e^{-\nu(v)(t-\tau)}\nu(v)|w_{{\ell_0} }\Pp f|^2\rmd v \rmd x\rmd \tau \notag\\
        \lesssim & \|w_{{\ell_0} }\Pp f_0\|^2+\int_0^t \int_{\mathbb{R}_v^3\times\mathbb{R}_x^3} e^{-\nu(v)(t-\tau)} [\P(v\cdot\nabla_x f) -(v\cdot\nabla_x\P f)] w_{{\ell_0} }^2 \Pp f\rmd v \rmd x\rmd \tau \notag\\
        &+\int_0^t \int_{\mathbb{R}_v^3\times\mathbb{R}_x^3} e^{-\nu(v)(t-\tau)} K\Pp f  w_{{\ell_0} }^2 \Pp f\rmd v \rmd x\rmd\tau +\int_0^t \int_{\mathbb{R}_v^3\times\mathbb{R}_x^3}  e^{-\nu(v)(t-\tau)} \Gamma(f,f) w_{{\ell_0} }^2 \Pp f\rmd v \rmd x\rmd\tau. \notag
    \end{align}
    Since \(\mathbf Pf\) is a linear combination of collision-invariant basis functions multiplied by \(\mu^{1/2}\), its Gaussian decay absorbs any polynomial velocity weight. Thus we have
    \begin{align}
       &\int_0^t \int_{\mathbb{R}_v^3\times\mathbb{R}_x^3} e^{-\nu(v)(t-\tau)} [\P(v\cdot\nabla_x f) -(v\cdot\nabla_x\P f)] w_{{\ell_0} }^2 \Pp f\rmd v \rmd x\rmd \tau \nonumber\\
       \lesssim&C_\eta\int_0^t \int_{\mathbb{R}_v^3}e^{-\nu(v)(t-\tau)}\mu^\delta [|\nabla_xf|_{L^2_x}^2+\|\nabla_xf\|_{L^2_{x,v}}^2]\rmd v\rmd \tau\nonumber\\
       &+\eta \int_0^t \int_{\mathbb{R}_v^3\times\mathbb{R}_x^3} e^{-\nu(v)(t-\tau)}\nu(v)|w_{\ell_0} \Pp f|^2\rmd v \rmd x\rmd \tau \notag.
    \end{align}
    The first term on the right-hand side has the time decay rates as $\|\nabla_xf\|^2\sim (1+t)^{-\frac 52}$. In some sense, this term replaces the estimate on
    $\int_0^t e^{-\lambda(t^p-\tau^p)}  |\xi|^2|\P \hat f|^2_{L^2_v}\rmd\tau$, induced by the time-velocity splitting method.

    However, there are extra difficulties arising from the estimate of $K$ and $\Gamma$ due to the presence of the factor $e^{-\nu(v)(t-\tau)}$. To this end, we apply the following strategy:
 
    \emph{(i) A hierarchy of time-weighted energy estimates.}
    We begin by deriving a weighted energy estimate for the microscopic component \(\mathbf P^\perp f\). More precisely, we multiply the microscopic equation by \((1+t)^{5/2}w_{\ell_0}^2\mathbf P^\perp f\) and integrate in $\mathbb{R}_v^3\times\mathbb{R}_x^3$. This leads to estimates of the form
    \begin{align}\label{k-estimate}
        &\int_0^t \int_{\mathbb{R}_v^3\times\mathbb{R}_x^3} e^{-\nu(v)(t-\tau)}(1+\tau)^{\frac 52} K\Pp f  w_{{\ell_0} }^2 \Pp f\rmd v \rmd x\rmd\tau\nonumber\\
        \lesssim&\eta\sup_{0\leq \tau\leq t} \left\{ (1+\tau)^{\frac52}\|w_{\ell_0} \Pp f(\tau)\|^2 \right\}+\cdots,
    \end{align} 
    where $\eta>0$ is sufficiently small and the omitted terms are controlled by the other estimates established below.
    
    A key observation is that differentiating the time weight produces a lower-order contribution with time weight \((1+t)^{3/2}\), that is 
    \[
        \int_0^t\int_{\mathbb R^3_v\times\mathbb R^3_x}
        e^{-\nu(v)(t-\tau)}(1+\tau)^{\frac32}
        |w_{\ell_0}\mathbf P^\perp f|^2\rmd v\rmd x\rmd \tau.
    \]
    Since the collision frequency satisfies $\nu(v)\sim \langle v\rangle^\gamma$ with $-3<\gamma<0$, 
	this contribution can be controlled by the microscopic dissipation by introducing the stronger velocity weight \(w_{\ell_0+1/2}\), for which
	\[
    	w_{\ell_0+1/2}
    	\sim\nu^{-1/2}(v)w_{\ell_0}.
	\]
    Repeating this argument at successive levels yields the time-weighted hierarchy
	\[
	\begin{aligned}
		&(1+t)^{5/2}
		\|w_{\ell_0}\mathbf P^\perp f\|^2,\\
		&(1+t)^{3/2}
		\|w_{{\ell_0}+1/2}\mathbf P^\perp f\|^2,\\
		&(1+t)^{1/2}
		\|w_{{\ell_0}+1}\mathbf P^\perp f\|^2
        .
	\end{aligned}
	\]
	The lower-order contribution generated by differentiating the third time weight carries the factor \((1+t)^{-1/2}\), which can be controlled by the corresponding energy estimate without a time weight, using the stronger velocity weight \(w_{\ell_0+3/2}\).
    Thus, no further time-weighted level is needed. 
    In this way, the hierarchy compensates for the loss generated by differentiating the time weights and allows us to retain the desired decay rate.
    
    \emph{(ii) A dyadic decomposition of the velocity space.}
    Indeed, due to the presence of the factor $e^{-\nu(v)(t-\tau)}$,  
    the preceding hierarchy does not suffice to control and absorb the contributions of \(K\) and \(\Gamma\) for the whole range of soft potentials. To obtain the desired estimate as \eqref{k-estimate}, we decompose velocity space into dyadic regions
    \[
        \mathbb{R}^3_v=\mathbb{V}_0 \cup \Big(\bigcup_{j\ge1}\mathbb V_j\Big),
    \]
    where \(\mathbb V_0\) is a bounded-velocity region and each \(\mathbb V_j\), \(j\geq1\), is an annulus. This localization allows us to exploit the velocity dependence of \(\nu(v)\) and apply the appropriate operator estimates on each region; see also \cite{Lei-Liu-Zhang-Zhao-2026-arxiv}.
    On \(\mathbb V_0\), the collision frequency is comparable to a positive constant, and the estimates for \(K\) and \(\Gamma\) follow from Lemmas \ref{lemma:K-R3} and \ref{lemma:Gamma-R3}. On each annulus \(\mathbb V_j\), we use the relation \(\nu(v)\sim\langle v\rangle^\gamma\) together with the uniform bound
	\[
	\int_0^t
	\mathbf 1_{\{v\in\mathbb V_j\}}
	e^{-\nu(v)(t-\tau)}\nu(v)\,\mathrm d\tau
	\leq C,
	\]
	and invoke the localized estimates in Lemmas \ref{lemma:K-Vj} and \ref{lemma:Gamma-Vj}. These estimates yield local bounds involving the following levels of the time-weighted hierarchy
	\begin{align*}
		&(1+t)^{5/2}
		\|w_{\ell_0}\mathbf P^\perp f\|_{L_v^2(\mathbb V_j)L_x^2(\mathbb{R}^3)},\\
		&(1+t)^{3/2}
		\|w_{\ell_0+1/2}\mathbf P^\perp f\|_{L_v^2(\mathbb V_j)L_x^2(\mathbb{R}^3)},\\
		&(1+t)^{1/2}
		\|w_{\ell_0+1}\mathbf P^\perp f\|_{L_v^2(\mathbb V_j)L_x^2(\mathbb{R}^3)}
        .
	\end{align*}
    The corresponding lower-order contribution generated by differentiating the third time weight is controlled analogously by the localized energy estimate without a time weight, with the velocity weight \(w_{\ell_0+3/2}\).
	Summing these localized estimates over the dyadic regions gives bounds compatible with the time-weighted hierarchy and makes it possible to absorb the relevant contributions of \(K\) and \(\Gamma\). 
    
    Combining these two ingredients with the decay estimates for the macroscopic component, we close the weighted estimates for the microscopic component. In particular, we obtain
    \[
        \|\mathbf P^\perp f\|^2
        \lesssim (1+t)^{-5/2} \epsilon_0^2,    \]
    which is exactly the optimal time decay rate of the microscopic part \(\mathbf P^\perp f\).

    \subsection{Notation}
        Before proceeding, we introduce some basic notation used throughout the paper:
        \begin{itemize}
        	\item $C$ denotes some positive (generally large) constant, and $\varepsilon$, $\eta$, $\lambda$, $\kappa_i$ and $\eta_i$ stand for some positive (generally small) constants. Note that all these constants may take different values in different places.
            
        	\item $A\lesssim B$ means that there is a generic constant $C>0$ such that $A\leq CB$. $A\sim B$ means $A\lesssim B$ and $B\lesssim A$.
            
        	\item For multi-indices $\alpha= [\alpha_1,\alpha_2, \alpha_3]$ and $\beta = [\beta_1, \beta_2, \beta_3]$, $\partial^{\alpha}_x\partial^{\beta}_v=\partial^{\alpha_1}_{x_1}\partial^{\alpha_2}_{x_2}\partial^{\alpha_3}_{x_3} \partial^{\beta_1}_{ v_1}\partial^{\beta_2}_{ v_2}\partial^{\beta_3}_{ v_3}$.
            The length of $\alpha$ is denoted by $|\alpha|=\alpha_1 +\alpha_2 +\alpha_3$. 
            $\alpha'\leq  \alpha$ means that no component of $\alpha'$ is greater than the corresponding component of $\alpha$, and $\alpha'<\alpha$ means that $\alpha'\leq  \alpha$ and $|\alpha'|<|\alpha|$.
            For an integer $k\ge0$, $\nabla_x^k f$ and $\nabla_v^k f$ denote the sums over all indices $\alpha$ and $\beta$ with $|\alpha|=k$ and $|\beta|=k$, respectively.
            
            \item $e_i$ denotes the multi-index with the $i$-th element unit and all other components equal to zero.

            \item For simplicity, we write \(L_v^2\) and \(L_x^2\) for \(L^2(\mathbb{R}_v^3)\) and \(L^2(\mathbb{R}_x^3)\), respectively. The notation \(\langle\cdot,\cdot\rangle\) denotes the \(L_v^2\) inner product, with the corresponding norm denoted by \(|\cdot|_{L_v^2}\), while \((\cdot,\cdot)\) denotes the \(L_{x,v}^2\) inner product, with the corresponding norm denoted by \(\|\cdot\|\).

        	\item For $1\leq p,q\leq \infty$, define the mixed velocity-space Lebesgue space $L_x^pL_v^q=L^p(\mathbb{R}_x^3;L^q(\mathbb{R}_v^3))$ with the norm
        	\begin{align*}
                \|g\|_{L^p_xL^q_v}=\left( \int_{\mathbb{R}^3_x} \left(\int_{\mathbb{R}^3_v}|g(x,v)|^q \rmd v\right)^{\frac pq} \rmd x \right)^{\frac 1p}.
        	\end{align*}
        
            \item  Define the following velocity-weighted norms involving the collision frequency:
            \begin{align*}
                \|g\|^2=\|g\|_{L^2_{x,v}}^2,
                \quad
                |g|_{\nu}^2=|\nu^\frac12g|_{L_v^2}^2,
                \quad 
                \|g\|_{\nu}^2=\|\nu^\frac12g\|^2.
            \end{align*}

            \item For $q\geq1$, define the mixed Lebesgue space $Z_q=L^2(\mathbb{R}^3_v;L^q(\mathbb{R}^3_x))$ with the norm
            \begin{align}
                \| g \|_{Z_q}=\left(\int_{\mathbb{R}_v^3} \left( \int_{\mathbb{R}_x^3}|g(x,v)|^q \mathrm{d}x \right)^\frac 2q \mathrm{d}v \right)^{\frac12}. \notag
            \end{align}
            
            \item Let $\mathscr{S}$ be the set of all Schwartz functions on $\mathbb{R}^3_x$. For any $g(x)\in\mathscr{S}(\mathbb{R}^3_x)$, its Fourier transform is defined by
            \begin{align*}
                \hat{g}(\xi) = \mathcal{F}g(\xi) = \int_{\mathbb{R}^3_x} e^{-\mathrm{i}x\cdot\xi}g(x) \mathrm{d}x, \quad x\cdot\xi = \sum_{j = 1}^3 x_j\xi_j,
            \end{align*}
            where $\mathrm{i} = \sqrt{-1} \in \mathbb{C}$.
        \end{itemize}
    \subsection{Outline of the paper}
    The remainder of the paper is organized as follows. Section~\ref{Sect-Pre} collects the linear and nonlinear estimates used in the analysis. Section~\ref{Sec:priori} establishes the a priori energy estimates and the basic decay bounds. Section~\ref{Sec:decayPp} proves the optimal temporal decay estimate for the microscopic component. For brevity, we defer some standard proof details to the appendix.

\section{Preliminaries}\label{Sect-Pre}
    In this section, we collect several estimates for the linearized and nonlinear collision operators, which will be used repeatedly in the proof of the main result.
    \subsection{Linear analysis}
    We begin by recalling the coercivity estimates for the linearized collision operator $L$.
    \begin{lemma}[\cite{Strain-Guo-2008-ARMA}]\label{lemma L estimate}
        (i) It holds that
        \begin{align} 
            \langle Lg, g \rangle \gtrsim  |\mathbf{P}^\perp g  |^2_{L^2_\nu}. \notag
        \end{align}
        (ii) Let $\ell\in\mathbb{R}$. There exists $C>0$ such that
        \begin{align} 
            \langle w_{\ell}^2Lg, g \rangle 
            \gtrsim 
            |w_\ell g |^2_{L^2_\nu}-C |\bar{\chi}_{C}g |^2_{L^2_v}.\notag
        \end{align}
        (iii) Let $\beta$ be a multi-index with $|\beta|>0$, and let $\ell\in\mathbb{R}$. For any $\eta>0$, there exists $C_\eta>0$ such that
        \begin{align}
            \langle w_{\ell-|\beta|}^2\partial^\beta_v [Lg], \partial^\beta_vg \rangle 
            \gtrsim 
             | w_{\ell-|\beta|} \partial^\beta_v g |^2_{L^2_\nu}
            -\eta\sum_{|\beta'| \leq |\beta|} |w_{\ell-|\beta'|} \partial^{\beta'}_v g |^2_{L^2_\nu}
            -C_\eta |\bar{\chi}_{C_\eta}(v)g |^2_{L^2_v}.\notag
        \end{align}
        Here $0\le\bar{\chi}_m(v)\le 1$ is a smooth cut-off function such that for $m\ge0$,
        \begin{align}
            \bar{\chi}_m(v)\equiv1,\quad |v|\le m;
            \quad 
            \bar{\chi}_m(v)\equiv0,\quad |v|\ge 2m.\notag
        \end{align}
    \end{lemma}

    For the operator $K$ on the whole velocity space, we use the following weighted estimate.
    \begin{lemma}[\cite{Strain-Guo-2008-ARMA}]\label{lemma:K-R3}
        Let $-3<\gamma<0$, $N\geq4$ and $\ell\geq N$. Let $0\leq |\beta|\leq N$. Then for all $\eta>0$, there exists $C_\eta>0$ such that
        \begin{align*}
            |\langle w^2_{\ell-|\beta|}\partial_v^{\beta}[Kg_1],g_2 \rangle|
            &\lesssim
            \left\{
            \eta\sum_{|\beta'|\le|\beta|} |w_{\ell-|\beta'|}\partial_v^{\beta'}g_1|_{L_\nu^2}
            +C_{\eta}|w_{\ell}\bar{\chi}_{C_\eta}(v)g_1|_{L_v^2}
            \right\}
            |w_{\ell}g_2|_{L_\nu^2}.
        \end{align*}
    \end{lemma}
    
    To estimate the contribution of \(K\) on each velocity annulus, we use the following localized estimate.
    \begin{lemma}[\cite{Lei-Liu-Zhang-Zhao-2026-arxiv}]\label{lemma:K-Vj}
        Let $-3<\gamma<0$, $N\geq4$ and $\ell\geq N$. Let $0\leq |\beta|\leq N$. Then it holds that
        \begin{align}
            \left|\int_{\mathbb{V}_j} \nu^{-1}(v)w_{\ell-|\beta|}^2\partial_v^\beta[K g_1]  g_2\mathrm{d}v\right|
            \lesssim& 
            \frac{\eta_1}{2^{2j}} \sum_{|\beta'|\le|\beta|} |{w}_{\ell-|\beta'|}\partial_v^{\beta'}g_1|_{L^2_v(\mathbb{R}^3)}^2
            +\frac{1}{2^{2j}} |g_1|_{L^2_v(\mathbb{R}^3)}^2 \notag\\
            &+\eta_2 \sum_{|\beta'|\le|\beta|} |w_{\ell-|\beta'|}\partial_v^{\beta'} g_1|_{L^2_v(\widetilde{\mathbb{V}}_j)}^2
            +\eta_3  |{w}_{\ell}g_2|_{L^2_v(\mathbb{V}_j)}^2,\notag
        \end{align}
        where $\widetilde{\mathbb{V}}_j=\mathbb{V}_{j-1}\cup \mathbb{V}_j\cup \mathbb{V}_{j+1}$ for $j\geq 1$, 
        and $\eta_i>0(i=1,2,3)$ are small constants.
    \end{lemma}
    
    \subsection{Nonlinear analysis}
    We now turn to estimates for the nonlinear collision operator $\Gamma$.
    Recall the definition of $\Gamma$ in \eqref{Gamma-def}. 
    As in \cite{Strain-Guo-2008-ARMA}, we have
    \begin{align}
        \partial^\alpha_x \partial^\beta_v \Gamma(g_1, g_2)
        \equiv& \sum C_\beta^{\beta_0 \beta_1 \beta_2} C_\alpha^{\alpha_1 \alpha_2}
        \Gamma^0  ( \partial_x^{\alpha_1} \partial_v^{\beta_1}g_1, \partial^{\alpha_2}_x \partial^{\beta_2}_v g_2 ) \notag\\
        \equiv& \int_{\mathbb{R}^3 \times \mathbb{S}^2} |v - u|^\gamma b(\cos\theta)
        \partial_v^{\beta_0} [ \mu^{\frac{1}{2}}(u)]\partial_x^{\alpha_1} \partial_v^{\beta_1}g_1(u') \partial^{\alpha_2}_x\partial^{\beta_2}_v g_2(v')\rmd u \rmd\omega\notag \\
        &-\partial^{\alpha_2}_x\partial^{\beta_2}_v g_2(v) \int_{\mathbb{R}^3 \times \mathbb{S}^2}|v - u|^\gamma b(\cos\theta) \partial_v^{\beta_0} [ \mu^{\frac{1}{2}}(u)]\partial_x^{\alpha_1} \partial_v^{\beta_1}g_1(u) \rmd u\rmd\omega, \notag
    \end{align}
    where the summation is over $\alpha_1 + \alpha_2 = \alpha$ and $\beta_0 + \beta_1 + \beta_2 = \beta$.
    We first recall the weighted estimates for the nonlinear collision operator \(\Gamma\) on the whole velocity space.
    \begin{lemma}[\cite{Guo-2003-ARMA,Strain-Guo-2008-ARMA, Duan-Yang-Zhao-2013-M3AS,Duan-Lei-Yang-Zhao-2017-CMP}]
        \label{lemma:Gamma-R3}
        Let $-3<\gamma<0$, $N\geq4$ and $\ell\geq N$. Let $\alpha_1+\alpha_2=\alpha$ and $\beta_0+\beta_1+\beta_2 =\beta$ with $0\leq|\alpha|+|\beta|\leq N$. Then it holds that
        \begin{align*}
            &|\langle w^2_{\ell-|\beta|} \Gamma^0(\partial_x^{\alpha_1}\partial_v^{\beta_1}g_1,\partial_x^{\alpha_2}\partial_v^{\beta_2}g_2),\partial_x^\alpha\partial_v^\beta g_3 \rangle| \\
            \lesssim& 
            \sum_{m\le2}\left\{ |\nabla_v^m (\mu^\delta\partial_x^{\alpha_1}\partial_v^{\beta_1}g_1 ) |_{L_v^2}
            + |w_{\ell-|\beta_1|} \partial_x^{\alpha_1}\partial_v^{\beta_1}g_1|_{L_v^2}  \right\}
            |w_{\ell-|\beta_2|} \partial_x^{\alpha_2}\partial_v^{\beta_2}g_2 |_{L_\nu^2}
            |w_{\ell-|\beta|} \partial_x^{\alpha}\partial_v^{\beta}g_3 |_{L_\nu^2},
        \end{align*}
        or
        \begin{align*}
            &|\langle w^2_{\ell-|\beta|} \Gamma^0(\partial_x^{\alpha_1}\partial_v^{\beta_1}g_1,\partial_x^{\alpha_2}\partial_v^{\beta_2}g_2),\partial_x^\alpha\partial_v^\beta g_3 \rangle| \\
            \lesssim& 
            \sum_{m\le2}\left\{ |\nabla_v^m (\mu^\delta\partial_x^{\alpha_2}\partial_v^{\beta_2}g_2 ) |_{L_v^2}+ |w_{\ell-|\beta_2|}\partial_x^{\alpha_2}\partial_v^{\beta_2}g_2 |_{L_v^2}\right\}
            |w_{\ell-|\beta_1|} \partial_x^{\alpha_1}\partial_v^{\beta_1}g_1 |_{L_\nu^2}
            |w_{\ell-|\beta|} \partial_x^{\alpha}\partial_v^{\beta}g_3 |_{L_\nu^2}.\notag
        \end{align*}
    \end{lemma}
    
    We next recall the localized estimate for \(\Gamma\) on each dyadic velocity annulus, which will be used in the energy estimates.
    \begin{lemma}[\cite{Lei-Liu-Zhang-Zhao-2026-arxiv}]
        \label{lemma:Gamma-Vj}
        Let $-3<\gamma<0$, $N\geq4$ and $\ell\geq N$. Let $\alpha_1+\alpha_2=\alpha$ and $\beta_0+\beta_1+\beta_2 =\beta$ with $0\leq|\alpha|+|\beta|\leq N$. Then for $j\in\mathbb{Z}^+$, it holds that
        \begin{align}
            &\left|\int_{\mathbb{V}_j}\nu^{-1}(v){w}_{\ell-|\beta|}^2\Gamma^0(\partial_x^{\alpha_1}\partial_v^{\beta_1}g_1,\partial_x^{\alpha_2}\partial_v^{\beta_2}g_2)\partial_x^{\alpha}\partial_v^{\beta}g_3\,\mathrm{d}v\right|\notag\\
            \lesssim& \sum_{m\le2} 
            \left\{  
            | \nabla_v^m (\mu^\delta\partial_x^{\alpha_1}\partial_v^{\beta_1}g_1) |_{L^2_v(\mathbb{R}^3)}
            +|{w}_{\ell-|\beta_1|}\partial_x^{\alpha_1}\partial_v^{\beta_1}g_1|_{L^2_v(\mathbb{R}^3)} 
            \right\} \notag\\
            &\times
            |{w}_{\ell-|\beta_2|}\partial_x^{\alpha_2}\partial_v^{\beta_2}g_2 |_{L^2_v(\widehat{\mathbb{V}}_j)}
            |{w}_{\ell-|\beta|}\partial_x^{\alpha}\partial_v^{\beta}g_3 |_{L^2_v(\mathbb{V}_j)},\notag
        \end{align}
        or
        \begin{align}
            &\left|\int_{\mathbb{V}_j}\nu^{-1}(v){w}_{\ell-|\beta|}^2\Gamma^0(\partial_x^{\alpha_1}\partial_v^{\beta_1}g_1,\partial_x^{\alpha_2}\partial_v^{\beta_2}g_2)\partial_x^{\alpha}\partial_v^{\beta}g_3\,\mathrm{d}v\right|\notag\\
            \lesssim& 
            \sum_{m\le2} 
            \left\{  
            | \nabla_v^m (\mu^\delta\partial_x^{\alpha_2}\partial_v^{\beta_2}g_2) |_{L^2_v(\mathbb{R}^3)}
            +|{w}_{\ell-|\beta_2|}\partial_x^{\alpha_2}\partial_v^{\beta_2}g_2|_{L^2_v(\mathbb{R}^3)} 
            \right\} \notag\\
            &\times
            |{w}_{\ell-|\beta_1|}\partial_x^{\alpha_1}\partial_v^{\beta_1}g_1 |_{L^2_v(\widehat{\mathbb{V}}_j)}
            |{w}_{\ell-|\beta|}\partial_x^{\alpha}\partial_v^{\beta}g_3 |_{L^2_v(\mathbb{V}_j)},\notag
        \end{align}
        where $\widehat{\mathbb{V}}_j=\mathbb{V}_{j-2}\cup \mathbb{V}_{j-1}\cup \mathbb{V}_j\cup \mathbb{V}_{j+1}$ for $j\geq 2$ 
        and $\widehat{\mathbb{V}}_1=\mathbb{V}_0\cup \mathbb{V}_{1}\cup \mathbb{V}_{2}$.
    \end{lemma}
    
    Finally, we record a basic $L_v^2$ estimate for $\Gamma$.
    \begin{lemma}[\cite{Duan-Yang-Zhao-2013-M3AS,Xiao-Xiong-Zhao-2013-JDE,Xiao-Xiong-Zhao-2017-JFA}]
        \label{lemma:Gamma-L2}
        Let $-3<\gamma<0$, $N\geq4$ and $\ell\geq N$. Then it holds that
        \begin{align*}
            | w_\ell\Gamma(g_1,g_2) |_{L_v^2} ^2
            &\lesssim 
            \sum_{m\le 2} 
            |w_{\ell-m}\nabla_v^m g_1 |_{L_\nu^2} ^2
            |w_\ell g_2 |_{L_\nu^2}^2,
        \end{align*}
        or
        \begin{align*}
            | w_\ell\Gamma(g_1,g_2) |_{L_v^2} ^2
            &\lesssim 
            \sum_{m\le 2} 
            |w_{\ell-m}\nabla_v^m g_2 |_{L_\nu^2}^2 
            |w_\ell g_1 |_{L_\nu^2}^2.
        \end{align*}
    \end{lemma}

\section{Global Well-posedness and basic time decay}\label{Sec:priori}
    In this section, we intend to deduce the uniform-in-time a priori estimates on the solution to the Cauchy problem \eqref{f-perturbation} and \eqref{f0}. 
    Let $-3 < \gamma < 0$, $N \geq 4$ and $\ell_0 \geq N$.
    Define the time-weighted energy norm $X(t)$ by
    \begin{align}\label{Xt}
        X(t) 
        &= \sup_{0 \leq \tau \leq t} \mathcal{E}_{N,\ell_0+1}(\tau) 
        + \sup_{0 \leq \tau \leq t} \left\{ (1+\tau)^{\frac{3}{2}} \mathcal{E}_{N,\ell_{0}}(\tau) \right\}.
    \end{align} 
    Suppose that the Cauchy problem  admits a smooth solution $f(t, x, v)$ over $0 \leq t \leq T$ for $0 < T \leq \infty$, and also the solution satisfies
    \begin{align}\label{priori assumption}
        \sup_{0 \leq t \leq T} X(t) \leq \delta_0,
    \end{align}
    where $\delta_0 > 0$ is a suitably small constant.

    \subsection{Basic time decay}
    Before carrying out the estimate, we apply the macro-micro decomposition. As in \cite{Duan-2008-JDE}, for $f(t,x,v)$, define moment functions $A_{mj}(f)$ and $B_j(f)$, $1\leq m, j\leq3$,  
    \begin{equation*}
        A_{mj}(f)=\int_{\mathbb{ R}^3}( v_m v_j-1)\mu^\frac12 f \rmd v,\quad B_j(f)=\frac{1}{10}\int_{\mathbb{ R}^3}(| v|^2-5) v_j\mu^\frac12 f \rmd v.
    \end{equation*}
    Then by taking the velocity integrations of \eqref{f-perturbation} with respect to the above high-order moments, we can derive a fluid-type system of equations
    \begin{equation}\label{Macro-equation}
        \begin{cases}
        \partial_ta+\nabla_x\cdot b=0,\\
        \partial_tb+\nabla_x(a+2c)+\nabla_x \cdot A({\bf P}^\perp f)=0,\\
        \partial_tc+\frac13\nabla_x\cdot b+\frac53\nabla_x\cdot B({\bf P}^\perp f)=0,
        \end{cases}
    \end{equation}
    and
    \begin{equation}\label{moment-equation}
        \begin{cases}
        \partial_tA_{mj}({\bf P}^\perp f)+\partial_{x_m}b_j+\partial_{x_j}b_m
        -\frac23\delta_{mj}\nabla_x\cdot B({\bf P}^\perp f)
        =A_{mj}(r+G),\\
        \partial_t B_{j}({\bf P}^\perp f)+\partial_{x_j}c=B_j(r+G),
        \end{cases}
    \end{equation}
    with
    \begin{equation}
        r=- v\cdot\nabla_x{\bf P}^\perp f-L{\bf P}^\perp f,\quad G={\Gamma}(f,f),\notag
    \end{equation}
    where $r$ is a linear term related only to the microscopic component $\mathbf{P}^\perp f$, and $G$ is a quadratic nonlinear term.

    The following lemma establishes the weighted energy estimate.
    \begin{lemma}\label{lemma:EE}
        For any $\ell\geq N$, there exists an energy functional $\mathcal{E}_{N,\ell}(t)$ satisfying \eqref{energy functional} such that
        \begin{align}\label{energy estimate}
            \frac{\rmd}{\rmd t}\mathcal{E}_{N,\ell} (t)+\lambda\mathcal{D}_{N,\ell}(t) \leq 0,
        \end{align}
        for any $0\leq t \leq T$ and some constant $\lambda>0$, where $\mathcal{D}_{N,\ell}(t)$ is defined in \eqref{dissipation functional}.
    \end{lemma}
    \begin{proof}
        Applying arguments similar to those in \cite{Liu-Yang-Yu-2004-PD, Guo-2003-ARMA}, we obtain the desired energy estimate. To keep the main text concise, we present the complete proof in the appendix.
    \end{proof}

    Now we consider the Cauchy problem on the linearized system with a non-homogeneous source term $h$:
    \begin{equation}\label{nonlinear system}
    	\left\{
    	\begin{aligned}
    		&\partial_t f+v \cdot \nabla_x f+Lf=h,\\
    		&f|_{t=0}=f_0,
    	\end{aligned}
    	\right.
    \end{equation} 
    where $h$ takes the form of
    \begin{align}
        h=\Gamma(f,f),\notag
    \end{align}
    which satisfies the mass conservation law 
    \[  
        \langle h,\sqrt{\mu} \rangle=0.
    \]

    For the linearized homogeneous Boltzmann equation, i.e. $h=0$, the solution to the Cauchy problem \eqref{nonlinear system} can be written as the mild form
    \begin{align}
        f=S(t)f_0,\notag
    \end{align}
    where $S(t)$ denotes the solution operator for the Cauchy problem on the linearized equation without any source. For the integer $k\geq0$, define the time decay rate index $\sigma_k$ by
    \begin{align}
        \sigma_k=\frac34+\frac k2. \notag
    \end{align}

    \begin{lemma}\label{E int}
        Let $f$ be the solution to the Cauchy problem \eqref{nonlinear system} with $h=0$. Then there exists a time-frequency interactive functional $\mathcal{E}_{int}(t,\xi)$ satisfying 
        \begin{align*}
            |\mathcal{E}_{int}(t,\xi)| \sim  |\hat{f}(\xi) |^2_{L^2_v},
        \end{align*}
        such that for any $t\geq 0$ and $\xi\in\mathbb{R}^3$, and for some $\lambda>0$,
        \begin{align}\label{interactive functional}
            \partial_t \mathcal{E}_{int}(t,\xi)
            + \frac{\lambda|\xi|^2}{1+|\xi|^2}   |(\widehat{a,b,c})|^2
            + \lambda|\mathbf{P}^\perp \hat{f}|_{L^2_\nu}^2
    		\leq 0,
        \end{align}
    \end{lemma}
    \begin{proof}
        Applying the Fourier transform to the first equation of $\eqref{nonlinear system}$ with $h=0$, one has 
        \begin{align}\label{f fourier}
            \partial_t \hat{f} + \mathrm{i} \xi \cdot v \hat{f} + L \hat{f}=0.
        \end{align}
        By taking the complex inner product of \eqref{f fourier} with $\hat{f}$, taking the real part and integrating it over $v\in\mathbb{R}^3$, one has 
        \begin{align}\label{f fourier homogeneous}
            \frac{1}{2}\partial_t |\hat{f}|^2_{L^2_v}
            +\lambda |\mathbf{P}^\perp \hat f|^2_{L^2_\nu} 
            \leq 0,
        \end{align}
        where we have used Lemma \ref{lemma L estimate}.

        By taking the Fourier estimate of the fluid-type equations \eqref{Macro-equation} and \eqref{moment-equation} with $G=0$, we can get the corresponding equations in $(t,\xi)$ and deduce as in \cite{Duan-2011-Nonlinearity} that 
        there exists an interactive functional $\mathcal{I}_{int}(t,\xi)$ satisfying 
        \begin{align*}
            |\mathcal{I}_{int}(t,\xi)| \lesssim  |\hat{f}(\xi) |^2_{L^2_v},
        \end{align*}
        such that
        \begin{align}\label{functional homogeneous}
            \partial_t \mathcal I_{\mathrm{int}}(t,\xi)
            + \frac{\lambda|\xi|^2}{1+|\xi|^2}  |\widehat{(a,b,c)}|^2 
    		\lesssim  
            |\mathbf{P}^\perp \hat{f}|_{L^2_\nu}^2.
        \end{align}
        Making a suitable combination of \eqref{f fourier homogeneous} and \eqref{functional homogeneous}, 
        we obtain the desired estimate \eqref{interactive functional},
        where $\mathcal{E}_{int}(t,\xi)$ is given by
        \begin{align*}
            \mathcal{E}_{int}(t,\xi) = |\hat{f}(\xi) |^2_{L^2_v}+\kappa_1\mathcal I_{\mathrm{int}}(t,\xi)\sim  |\hat{f}(\xi) |^2_{L^2_v}, 
        \end{align*}
        for a constant $\kappa_1>0$ small enough.
        This completes the proof of Lemma \ref{E int}.
    \end{proof}

    To obtain the optimal time decay estimates, we first establish decay bounds for the solution and its first-order spatial derivatives.
    \begin{lemma}\label{lemma:nablaf}
        Let $-3<\gamma<0$ and $N\geq4$. 
        Let $f(t,x,v)$ be the solution to the Cauchy problem \eqref{f-perturbation} and \eqref{f0}. For $k=0,1$, it holds that
        \begin{align}
			\| \nabla_x^k f(t) \|
			\lesssim (1+t)^{-\frac34-\frac{k}{2}}
			\left\{ \| w_{\frac{N}{2}} \nabla_x^k f_0 \|+\| w_{\frac{N}{2}} f_0 \|_{Z_1}+\sup_{t}\left\{(1+t)^\frac{3}{2}\mathcal{E}_{N,N}\left(t\right)\right\}\right\}.\notag
		\end{align}
    \end{lemma}
    \begin{proof}
        We first consider the linear problem. 
        By taking the complex inner product of \eqref{f fourier} with $w_{\ell}^2 \hat{f}$, taking the real part and integrating it over $v\in\mathbb{R}^3$, one has for some $\lambda >0$,
        \begin{align}\label{wf fourier homogeneous}
            \frac{1}{2}\partial_t  |w_\ell \hat{f} |^2_{L^2_v}
            +  \lambda|w_\ell \hat f|^2_{L^2_\nu} 
            \lesssim 
            |\hat f|^2_{L^2_\nu},
        \end{align}
        where we have used Lemma \ref{lemma L estimate} and 
        $|\bar{\chi}_Cw_\ell\hat{f}|^2_{L^2_v} \lesssim |\hat f|^2_{L^2_\nu}$.

        Applying $\mathbf{P}^\perp$ to \eqref{f fourier}, we have
        \begin{align}\label{P perp f fourier}
            \partial_t\mathbf{P}^{\perp} \hat{f} +\mathrm{i} \xi \cdot v \mathbf{P}^{\perp}\hat{f} +L\mathbf{P}^\perp \hat{f} = -\mathrm{i} \xi \cdot v \mathbf{P}\hat{f}+{\bf P}(\mathrm{i} \xi \cdot v \hat{f}).
        \end{align}
        By taking the complex inner product of \eqref{P perp f fourier} with $w_{\ell}^2\mathbf{P}^\perp \hat{f}$, taking the real part and integrating it over $v\in\mathbb{R}^3$, one has for some $\lambda >0$,
        \begin{align}\label{fourier weighted P perp f}
            &\frac{1}{2}\partial_t |w_\ell\mathbf{P}^\perp \hat{f} |^2_{L^2_v}
            +\lambda  |w_\ell \Pp \hat f|^2_{L^2_\nu}
            \notag\\
            \lesssim  
            &\mathrm{Re} \langle -\mathrm{i} \xi \cdot v \mathbf{P}\hat{f},w_{\ell}^2\Pp \hat{f}\rangle 
            +\mathrm{Re} \langle {\bf P}(\mathrm{i} \xi \cdot v \hat{f}),w_{\ell}^2\Pp \hat{f}\rangle 
            + |\Pp \hat f|^2_{L^2_\nu},
        \end{align}
        where we have used Lemma \ref{lemma L estimate} and 
        $| \bar{\chi}_Cw_\ell\Pp \hat{f} |^2_{L^2_v} \lesssim | \Pp \hat f|^2_{L^2_\nu}$.
        For the first term on the right-hand side of \eqref{fourier weighted P perp f}, one has
        \begin{align}
            | \mathrm{Re} \langle -\mathrm{i} \xi \cdot v \mathbf{P}\hat{f}, w_{\ell}^2\Pp \hat f \rangle |  \nonumber
            &\leq \int_{\mathbb{R}^3_v} |\xi| |v| |(\widehat{a,b,c})| \mu^{\delta} w_{\ell}^2 |\Pp \hat f|\rmd v \nonumber\\
            &\leq \int_{\mathbb{R}^3_v} |\xi|  |(\widehat{a,b,c})\langle v\rangle^\frac{\gamma}{2} |\mu^{\delta'}  |w_{\ell} \Pp \hat f \langle v \rangle^\frac{\gamma}{2} |\rmd v \nonumber\\
            &\leq 
            |\xi|^2 
            |\hat f |^2_{L^2_\nu}
            +\eta 
            | w_\ell \Pp \hat f  |^2_{L^2_\nu}.\notag
        \end{align}
        For the second term on the right-hand side of \eqref{fourier weighted P perp f}, recall the definition of ${\bf P}$ and denote 
		\begin{align}
			\phi_i=\mu^\delta, v_i\mu^\delta, (|v|^2-3)\mu^\delta, \quad i=1,2,3, \label{0phi-i}
		\end{align} 
        which implies
        \begin{align}
             |{\bf P}(\mathrm{i} \xi \cdot v \hat{f}) |
            \lesssim \int_{\mathbb{R}_v^3} |\xi||v|  |\hat{f} | |\phi_i| \mathrm{d}v\mu^\delta 
            \lesssim |\xi| 
            |\hat{f}|_{L^2_\nu}
            \mu^\delta.\notag
        \end{align}
        Then one has
        \begin{align}
            | \mathrm{Re} \langle {\bf P}(\mathrm{i} \xi \cdot v \hat{f}),w_{\ell}^2\Pp \hat f \rangle |
            \lesssim& \int_{\mathbb{R}^3_v} 
            |\xi| |\hat{f}|_{L^2_\nu}
            \mu^{\delta'} |w_\ell\Pp \hat f\langle v\rangle^\frac{\gamma}{2} | \rmd v \notag\\
            \lesssim& |\xi|^2 
            |\hat{f}|_{L^2_\nu}^2
            +\eta |w_\ell\Pp \hat{f}|_{L^2_\nu}^2. \notag
        \end{align}
        Plugging the above estimates into \eqref{fourier weighted P perp f}, one has
        \begin{align}\label{weighted P perp f estimate}
            \partial_t |w_\ell \mathbf{P}^\perp \hat{f} |^2_{L^2_v}
            +\lambda |w_\ell \mathbf{P}^\perp \hat f|^2_{L^2_\nu} 
            \lesssim 
            |\xi|^2 |\hat f |_{L^2_\nu}^2
            + |\Pp \hat f |^2_{L^2_\nu}.
        \end{align}
        
        In the case when $|\xi|\geq 1$ which implies $|\xi|^2/(1+|\xi|^2)\geq 1/2$, one can combine \eqref{interactive functional} and \eqref{wf fourier homogeneous} to obtain
        \begin{align}
            \partial_t \mathcal{E}_\ell^H(t,\xi)+  
            \lambda|w_\ell \hat{f} |^2_{L^2_\nu} 
            \leq 0, \notag
        \end{align}
        where $\mathcal{E}_\ell^H(t,\xi)$ is defined by
        \begin{align}
            \mathcal{E}_\ell^H(t,\xi)=\left\{ \mathcal{E}_{int}(t,\xi)+\kappa_2 |w_\ell \hat{f} |^2_{L^2_v}  \right\}\chi_{|\xi|\geq 1},\notag
        \end{align}
        for a constant $\kappa_2>0$ small enough.
        In the case when $|\xi|\leq 1$ which implies $|\xi|^2/(1+|\xi|^2)\geq |\xi|^2/2$, one can combine \eqref{interactive functional} and \eqref{weighted P perp f estimate} to obtain
        \begin{align}
            \partial_t \mathcal{E}_\ell^L(t,\xi)
            + \lambda|\xi|^2  
            |w_\ell \hat f|^2_{L^2_\nu}  
            \leq 0,\notag
        \end{align}
        where $\mathcal{E}_\ell^L(t,\xi)$ is defined by
        \begin{align}
            \mathcal{E}_\ell^L(t,\xi)
            =\left\{ \mathcal{E}_{int}(t,\xi)+\kappa_3 |w_\ell \mathbf{P}^\perp \hat{f} |^2_{L^2_v} \right\}\chi_{|\xi|\leq 1},\notag
        \end{align}
        for a constant $\kappa_3>0$ small enough. 
        Consequently, for $\ell\geq 0$, by introducing
		\begin{align}
		    \mathcal E_\ell(t,\xi) 
            =\mathcal E_\ell^H(t,\xi)+\mathcal E_\ell^L(t,\xi)
            \sim  |w_\ell \hat f |_{L^2_v}^2,\notag
		\end{align}
        we get the following pointwise in $\xi$ dissipative inequality
        \begin{align}
		      \partial_t\mathcal E_\ell(t,\xi)
		      + \frac{\lambda|\xi|^2}{1+|\xi|^2}  |w_{\ell}\hat f  |_{L^2_\nu}^2   \le 0, \notag
		\end{align}
        that is
		\begin{align}
		      \partial_t\mathcal E_\ell(t,\xi)
		      + \frac{\lambda|\xi|^2}{1+|\xi|^2}  |w_{\ell-\frac12}\hat f |_{L^2_v}^2   \le 0. \label{dissipation}
		\end{align}
		Note that $\mathcal E_\ell(t,\xi)$ is nonincreasing in time $t$.
        
        For $j>0$, the interpolation inequality gives that 
        \begin{align}
            \mathcal E_\ell(t,\xi) 
            \lesssim \mathcal E_{\ell-\frac12}^{\frac{j}{j+1}}(t,\xi)\mathcal E_{\ell+\frac{j}2}^{\frac{1}{j+1}}(t,\xi), \notag
        \end{align}
        which implies
		\[
		      \mathcal E_\ell^{\frac{j+1}{j}}(t,\xi) 
            \lesssim \mathcal E_{\ell-\frac12}(t,\xi) \mathcal E_{\ell+\frac{j}{2}}^{\frac{1}{j}}(t,\xi)
            \lesssim  |w_{\ell-\frac12}\hat f |^2_{L^2_v} \mathcal E_{\ell+\frac{j}{2}}^{\frac{1}{j}}(0,\xi). 
		\]
		Together with \eqref{dissipation}, one has
		\[
		    \partial_t \mathcal E_\ell(t,\xi)+ \frac{\lambda|\xi|^2}{1+|\xi|^2}\mathcal E_{\ell+\frac{j}2}^{-\frac{1}{j}}(0,\xi) \mathcal E_\ell^{\frac{j+1}{j}}(t,\xi)  \le 0.
		\]
		Solving the above inequality yields that
		\begin{align}
		    \mathcal E_\ell(t,\xi) 
            \lesssim  \mathcal E_{\ell+\frac j2}(0,\xi) \left( 1 + \frac{|\xi|^2}{1+|\xi|^2}\cdot \frac{\lambda t}{j} \right)^{-j} .\label{El 1}
		\end{align}
        
		Multiplying \eqref{El 1} by $|\xi|^{2k}$ and integrating it over $|\xi|\ge 1$, one has
        \begin{align}
            \int_{|\xi| \geq 1} |\xi|^{2k} \mathcal{E}_\ell(t, \xi) \rmd\xi
            \lesssim& \left(1 + \frac{\lambda t}{2j}\right)^{-j} \int_{|\xi| \geq 1} |\xi|^{2k} \mathcal{E}_{\ell+\frac{j}{2}}(0, \xi) \rmd\xi\notag\\
            \lesssim& (1+t)^{-j} \| w_{\ell+\frac{j}{2}} \nabla_x^k f_0 \|^2.\label{xi>1}
        \end{align}
        While for $|\xi|\leq 1$, one can take 
        $j>2\sigma_k=\frac32+k$ so that
        \begin{align}
            \int_{|\xi| \leq 1} |\xi|^{2k} \mathcal{E}_\ell(t, \xi) \rmd\xi
            \lesssim& |\mathcal{E}_{\ell+\frac{j}{2}}(0, \xi)|_{L^\infty_\xi} \int_{|\xi| \leq 1} |\xi|^{2k} \left(1 + \frac{\lambda|\xi|^2}{2j}t\right)^{-j}   \rmd\xi \notag\\
            \lesssim& (1+t)^{-\frac32-k} \| w_{\ell+\frac{j}{2}} f_0 \|_{Z_1}^2 .\label{xi<1}
        \end{align}
		By taking $j=2+k$, combining \eqref{xi>1} and \eqref{xi<1} yields that
		\begin{align*}
		    \| w_\ell \nabla_x^k S(t) f_0\|^2
    		&\lesssim 
            (1+t)^{-\frac32-k} \left( \| w_{\ell+\frac k2+1} \nabla_x^k f_0\|^2 + \| w_{\ell+\frac k2+1} f_0\|_{Z_1}^2 \right),
		\end{align*}
		which implies
        \begin{align}
            \| \nabla_x^k S(t) f_0\|^2
    		&\lesssim 
            (1+t)^{-\frac32-k} \left( \| w_{\frac k2+1} \nabla_x^k f_0\|^2 + \| w_{\frac k2+1} f_0\|_{Z_1}^2\right).\notag
        \end{align}

        Now we consider the nonlinear problem. By using the Duhamel principle, the solution to the Cauchy problem \eqref{nonlinear system} can be denoted as
        \begin{align*}
            f(t,x,v) = S(t)f_0 + \int_0^t S(t-\tau)h(\tau)\rmd\tau,
        \end{align*}
        with $h=\Gamma(f,f)$, which implies
        \begin{align}\label{decay111}
            \| \nabla_x^k f(t) \|
            \lesssim& (1+t)^{-\frac{3}{4}-\frac{k}{2}}  \left( \|w_{\frac k2+1}\nabla_x^kf_0\| + \|w_{\frac k2+1}f_0\|_{Z_1} \right) \notag\\
            &+ \int_0^t (1+t-\tau)^{-\frac{3}{4}-\frac{k}{2}} \left( \|w_{\frac k2+1}\nabla_x^k\Gamma(f,f)(\tau) \| + \|w_{\frac k2+1}\Gamma(f,f)(\tau) \|_{Z_1} \right) \rmd \tau. 
        \end{align} 
        Lemma \ref{lemma:Gamma-L2} yields that for $\alpha_1+\alpha_2=\alpha$ with $|\alpha|=k$, 
        \begin{align}
            \left| w_{\ell}\Gamma^0(\partial^{\alpha_1}_xf,\partial^{\alpha_2}_xf) \right|_{L^2_v}^2
            \lesssim 
            \sum_{m\le 2} 
             |w_{\ell-m}\nabla_v^m\partial^{\alpha_1}_x f |_{L_v^2}^2 
             |w_{\ell} \partial^{\alpha_2}_xf  |_{L_v^2}^2, 
            \label{Gg1}
        \end{align}
        or
        \begin{align}
            | w_{\ell}\Gamma^0(\partial^{\alpha_1}_xf,\partial^{\alpha_2}_xf) |_{L^2_v}^2
            \lesssim 
            \sum_{m\le 2}  | w_{\ell-m}\nabla_v^m\partial^{\alpha_2}_x f |_{L_v^2}^2 
            | w_{\ell} \partial^{\alpha_1}_xf  |_{L_v^2}^2. 
            \label{Gg2}
        \end{align}
        Combining the estimates above, we obtain
        \begin{align}
            \| w_{\frac k2+1}\Gamma(f,f) \|_{Z_1}
            \lesssim& 
            \sum_{m\le 2} 
            \| w_{\frac k2+1-m}\nabla_v^m f\|
            \| w_{\frac k2+1} f \|
            \lesssim
            \mathcal{E}_{N,N}(t). \notag
        \end{align}
        Moreover, for $0\leq |\alpha| \leq 1$, if $|\alpha_1|=0$, it follows from \eqref{Gg1} and the $L^\infty-L^2-L^2$ estimate that 
        \begin{align}
            \| w_{\frac k2+1}\Gamma^0(\partial^{\alpha_1}_x f, \partial^{\alpha_2}_x f) \|^2
            \lesssim& 
            \sum_{m\le 2} 
            \| w_{\frac k2+1-m}\nabla_v^m f\|_{L^\infty_xL^2_v}^2
            \| w_{\frac k2+1} \partial^{\alpha}_xf \|^2
            \notag\\
            \lesssim&
            \sum_{m,n\le 2} 
            \| w_{\frac k2+1-m}\nabla_x^n\nabla_v^m f\|^2
            \| w_{\frac k2+1} \partial^{\alpha}_xf \|^2\notag\\
            \lesssim&
            \mathcal{E}_{N,N}^2(t) . \notag
        \end{align}
        If $|\alpha_1|=|\alpha|$, the same estimate follows from \eqref{Gg2} and the $L^\infty-L^2-L^2$ estimate. 
        Consequently, for $k=0,1$, we obtain 
        \begin{align}\label{decay222}
            &\int_0^t (1+t-\tau)^{-\frac{3}{4}-\frac{k}{2}} \left( \|w_{\frac k2+1}\nabla_x^k\Gamma(f,f)(\tau) \| + \|w_{\frac{k}{2}+1}\Gamma(f,f)(\tau) \|_{Z_1}  \right) \rmd \tau\notag\\
            \lesssim& \int_0^t (1+t-\tau)^{-\frac{3}{4}-\frac{k}{2}} (1+\tau)^{-\frac32}   (1+\tau)^{\frac{3}{2}} \mathcal{E}_{N, N}(\tau)  \rmd \tau \notag\\
            \lesssim& (1+t)^{-\frac{3}{4}-\frac{k}{2}} \sup_{t} \left\{ (1+t)^{\frac{3}{2}} \mathcal{E}_{N,N}(t) \right\}.
        \end{align}
        Substituting \eqref{decay222} into \eqref{decay111}, for $k=0,1$, we obtain
        \begin{align}
            \| \nabla_x^k f(t) \|
            \lesssim& (1+t)^{-\frac{3}{4}-\frac{k}{2}}  \left( \|w_{\frac k2+1}\nabla_x^kf_0\| +\|w_{\frac k2+1}f_0\|_{Z_1}+ \sup_{t} \left\{ (1+t)^{\frac{3}{2}} \mathcal{E}_{N,N}(t) \right\}\right) \notag\\
            \lesssim& (1+t)^{-\frac{3}{4}-\frac{k}{2}}  \left( \|w_{\frac N2}\nabla_x^kf_0\| +\|w_{\frac N2}f_0\|_{Z_1}+ \sup_{t} \left\{ (1+t)^{\frac{3}{2}} \mathcal{E}_{N,N}(t) \right\}\right). \notag 
        \end{align}
        This completes the proof of Lemma \ref{lemma:nablaf}. 
    \end{proof}

    Based on Lemma \ref{lemma:nablaf}, we can further obtain the time decay of energy functional $\mathcal{E}_{N,\ell_0}(t)$.
    \begin{lemma}\label{lemma:k=0}
        Let $N\geq4$ and $\ell_0\geq N$. Let $T>0$. Let $f$ be the solution to the Cauchy problem \eqref{f-perturbation} and \eqref{f0} on $[0,T]$.
        Then for $0\leq t\leq T$, it holds that
        \begin{align}
            (1+t)^\frac32 \mathcal{E}_{N,\ell_0}(t) 
            \lesssim 
            \mathcal{E}_{N,\ell_0+1}(0) 
            + \| w_{\frac N2} f_0 \|^2_{Z_1} 
            + \sup_{0\leq\tau\leq t} \left\{ (1+\tau)^\frac32 \mathcal{E}_{N,N}(\tau)\right\}^2.\notag
        \end{align}
    \end{lemma}
    \begin{proof}
        Let $0<\delta<\frac12$. 
        Multiplying \eqref{energy estimate} with $\ell=\ell_0$ by $(1+t)^{\frac 32+\delta}$, one has for some $\lambda >0$,
        \begin{align}
            \frac{\rmd}{\rmd t} \left\{ (1+t)^{\frac 32+\delta}\mathcal{E}_{N,\ell_0}(t) \right\}
            +\lambda(1+t)^{\frac 32+\delta}\mathcal{D}_{N,\ell_0}(t) 
            \lesssim 
            (1+t)^{\frac 12+\delta}\mathcal{E}_{N,\ell_0}(t). \label{E1}
        \end{align}
        Recall \eqref{energy functional} and \eqref{dissipation functional}. Thus one has
        \begin{align}\label{E-D}
            \mathcal{E}_{N,\ell}(t) 
            \lesssim \mathcal{D}_{N,\ell+\frac12} (t) +\| \P f\|^2
            \lesssim \mathcal{D}_{N,\ell+\frac12} (t) +\| f\|^2. 
        \end{align}
        Combining \eqref{E1} and \eqref{E-D} with $\ell=\ell_0$, one has
        \begin{align}\label{EEE1}
            &\frac{\rmd}{\rmd t} \left\{ (1+t)^{\frac 32+\delta}\mathcal{E}_{N,\ell_0}(t) \right\}
            +\lambda(1+t)^{\frac 32+\delta}\mathcal{D}_{N,\ell_0}(t) \notag\\
            \lesssim &
            (1+t)^{\frac 12+\delta}\mathcal{D}_{N,\ell_0+\frac12}(t)
            +(1+t)^{\frac 12+\delta} \| f\|^2.
        \end{align}
        In a similar way, multiplying \eqref{energy estimate} with $\ell=\ell_0+\frac12$ by $(1+t)^{\frac 12+\delta}$ gives
        \begin{align}\label{EEE2}
            &\frac{\rmd}{\rmd t} \left\{ (1+t)^{\frac 12+\delta}\mathcal{E}_{N,\ell_0+\frac12}(t) \right\}
            +\lambda(1+t)^{\frac 12+\delta}\mathcal{D}_{N,\ell_0+\frac12}(t) 
            \notag\\
            \lesssim & (1+t)^{-\frac 12+\delta}\mathcal{D}_{N,\ell_0+1}(t)
            +(1+t)^{-\frac 12+\delta} \| f\|^2.
        \end{align}
        By taking $\ell=\ell_0+1$ in \eqref{energy estimate}, one has
        \begin{align}\label{EEE3}
            \frac{\rmd}{\rmd t} \mathcal{E}_{N,\ell_0+1}(t) 
            +\lambda \mathcal{D}_{N,\ell_0+1}(t) 
            \leq 0.  
        \end{align}
        Making a suitable combination of \eqref{EEE1}, \eqref{EEE2} and \eqref{EEE3}, one has
        \begin{align}
            &\frac{\rmd}{\rmd t} \left\{ (1+t)^{\frac 32+\delta}\mathcal{E}_{N,\ell_0}(t)+C_1(1+t)^{\frac 12+\delta}\mathcal{E}_{N,\ell_0+\frac12}(t) + C_2\mathcal{E}_{N,\ell_0+1}(t) \right\} \notag\\
            &+ \lambda \left\{ (1+t)^{\frac 32+\delta}\mathcal{D}_{N,\ell_0}(t)+C_1(1+t)^{\frac 12+\delta}\mathcal{D}_{N,\ell_0+\frac12}(t) +C_2\mathcal{D}_{N,\ell_0+1}(t) \right\}
            \notag\\
            \lesssim&
            (1+t)^{\frac12+\delta} \| f(t) \|^2, \notag
        \end{align}
        where $C_2\gg C_1 \gg 1$.
        Together with Lemma \ref{lemma:nablaf}, integrating the above inequality with respect to $t$ yields that
        \begin{align*}
            &(1+t)^{\frac 32+\delta}\mathcal{E}_{N,\ell_0}(t)
            \notag\\
            \lesssim&
            \mathcal{E}_{N,\ell_0+1}(0)
            +\int_0^t (1+\tau)^{-1+\delta} \left\{ \| w_{\frac N2} f_0 \|^2 + \| w_{\frac N2} f_0 \|^2_{Z_1} +  \sup_{0\leq\tau\leq t} \left\{ (1+\tau)^\frac32 \mathcal{E}_{N,N}(\tau)\right\}^2\right\} \rmd \tau \notag\\
            \lesssim&
            \mathcal{E}_{N,\ell_0+1}(0)
            +(1+t)^{\delta} \left\{ \| w_{\frac N2} f_0 \|^2 + \| w_{\frac N2} f_0 \|^2_{Z_1} +  \sup_{0\leq\tau\leq t} \left\{ (1+\tau)^\frac32 \mathcal{E}_{N,N}(\tau)\right\}^2\right\},
        \end{align*}
        which implies
        \begin{align*}
            (1+t)^{\frac 32}\mathcal{E}_{N,\ell_0}(t)
            \lesssim 
            \mathcal{E}_{N,\ell_0+1}(0) + \| w_{\frac N2} f_0 \|^2_{Z_1} + \sup_{0\leq\tau\leq t} \left\{ (1+\tau)^\frac32 \mathcal{E}_{N,N}(\tau)\right\}^2.
        \end{align*}
        This completes the proof of Lemma \ref{lemma:k=0}.
    \end{proof}

    In the following proposition, we  derive the a priori estimate for the Cauchy problem \eqref{f-perturbation} and \eqref{f0}.
    \begin{proposition}\label{proposition}
        Let $N\geq4$ and $f_0$, $\epsilon_0$ be defined as in Theorem \ref{thm}. Let $T>0$ and $f(t,x,v)$ be the solution to the Cauchy problem \eqref{f-perturbation} and \eqref{f0} on $[0,T]$. Recall \eqref{Xt}. If $X(t)$ is sufficiently small, then it holds that
        \begin{align}
            X(t)\lesssim \epsilon_0^2, \label{XfY}
        \end{align}
        for any $0 \leq t \leq T$.
    \end{proposition}
    \begin{proof}
        From Lemma \ref{lemma:EE}, one has 
		\begin{align}
			\frac{\rmd}{\rmd t} \mathcal{E}_{N,\ell_0+1}(t) + \lambda \mathcal{D}_{N,\ell_0+1}(t) \leq 0, \notag
		\end{align}
		which implies
		\begin{align}\label{E down}
			\mathcal{E}_{N,\ell_0+1}(t) \leq \mathcal{E}_{N,\ell_0+1}(0) \leq \epsilon_0^2,
		\end{align}
		for any $0\leq t \leq T$.
        From Lemma \ref{lemma:k=0}, one has
		\begin{align}\label{EN}
			(1+t)^{\frac 32}\mathcal{E}_{N,\ell_0}(t)
			\lesssim 
			\mathcal{E}_{N,\ell_0+1}(0) + \| w_{\frac N2} f_0 \|^2_{Z_1} + \sup_{0\leq\tau\leq t} \left\{ (1+\tau)^\frac32 \mathcal{E}_{N,\ell_0}(\tau)\right\}^2
			\lesssim
			\epsilon_0^2+X^2(t).
		\end{align}
		Recall \eqref{Xt}.
        Combining \eqref{E down} and \eqref{EN} gives that
        \begin{align*}
            X(t) \lesssim \epsilon_0^2 +X^2(t).
        \end{align*}
        It follows immediately that \eqref{XfY} holds for any $0\leq t \leq T$ for $\epsilon_0$ sufficiently small.
		This completes the proof of Proposition \ref{proposition}.
    \end{proof}

    \subsection{Global existence}
    Recall \eqref{Xt}. It holds that $X(0)\lesssim \epsilon_0^2$. By a standard argument, we have the following local existence result. If $X(0)$ is small enough, then there is a positive time $T^*>0$ such that \eqref{f-perturbation} and \eqref{f0} admit a solution in $[0,T^*]$ with $X(T^*)\leq 2X(0)$.

    From Proposition \ref{proposition}, if $\epsilon_0$ is sufficiently small, by the above local estimate and the usual continuity argument, the Cauchy problem \eqref{f-perturbation} and \eqref{f0} admit a unique global solution verifying
    \begin{align}
        \mathcal{E}_{N,\ell_0}(t) \lesssim (1+t)^{-\frac32}\epsilon_0^2. \notag
    \end{align} 

\section{Optimal time decay rate for microscopic part}\label{Sec:decayPp}
    In this section, we establish the optimal time decay estimate for the microscopic part $\mathbf{P}^{\perp}f$.
    
    \subsection{Microscopic part without weight}
    We first derive a time-decay estimate for the microscopic part $\mathbf{P}^{\perp}f$ without velocity weight. This estimate will serve as a preliminary step toward establishing the corresponding weighted decay estimates.
    \begin{lemma}\label{lemma Pp f}
        Let $-3<\gamma<0$ and $N\geq4$. Let $f(t,x,v)$ be the solution to  the Cauchy problem \eqref{f-perturbation} and \eqref{f0}. Then it holds that
        \begin{align}
            \sup_t \left\{ (1+t)^\frac52 \|\mathbf{P}^{\perp}f\|^2 \right\}
            &\lesssim \|\mathbf{P}^{\perp}f_0\|^2
            +\eta \sup_t \left\{ (1+t)^{\frac52} \|\mathbf{P}^{\perp}f \langle v\rangle^{-\frac{\gamma}{2}}\|^2 \right\}
            +\sup_t\left\{(1+t)^{\frac{5}{2}}\|\nabla_x f\|^2 \right\}, \notag
        \end{align}
        where $\eta>0$ is a small constant.
    \end{lemma}
    \begin{proof}
        Applying $\mathbf{P}^{\perp}$ to \eqref{f-perturbation}, we derive the microscopic equation
        \begin{align}\label{micro-eq}
            \partial_t\mathbf{P}^{\perp}f + v\cdot\nabla_x\mathbf{P}^\perp f + L\mathbf{P}^{\perp}f = -v\cdot\nabla_x\mathbf{P}f + {\bf P}(v\cdot\nabla_xf) +\Gamma(f,f).
        \end{align}
        Multiplying \eqref{micro-eq} by $\mathbf{P}^{\perp}f$ and integrating it over $\mathbb{R}_v^3\times\mathbb{R}_x^3$, we obtain
        for some $\lambda >0$,
        \begin{align}\label{Pp f energy}
            \frac{\mathrm{d}}{\mathrm{d}t}\|\mathbf{P}^{\perp}f\|^2
            + \lambda\|\mathbf{P}^{\perp}f\langle v\rangle^{\frac{\gamma}{2}}\|^2
            \lesssim & 
            - ( v\cdot\nabla_x\mathbf{P}f, \mathbf{P}^{\perp}f  )
            + ( {\bf P}(v\cdot\nabla_x f), \mathbf{P}^{\perp}f ) 
            + ( \Gamma(f,f),\mathbf{P}^{\perp}f  ), 
        \end{align}
        where we have used the fact that $(v\cdot\nabla_x\mathbf{P}^\perp f, \Pp f)=0$ and the coercivity estimate in Lemma \ref{lemma L estimate}.
        For the first term on the right-hand side of \eqref{Pp f energy}, using Young's inequality, we obtain
        \begin{align}
            ( v\cdot\nabla_x\mathbf{P}f, \mathbf{P}^{\perp}f  )
            \lesssim  ( \nabla_x(a,b,c) \mu^\delta, \mathbf{P}^{\perp}f  )
            \lesssim \|\nabla_x f\|^2 + \eta \|\Pp f\|^2. \notag
        \end{align}
        The second term vanishes by the orthogonality of $\mathbf{P}$ and $\mathbf{P}^{\perp}$, that is
        \begin{align}
            ( {\bf P}(v\cdot\nabla_x f), \mathbf{P}^{\perp}f )=0. \notag
        \end{align}
        To estimate the nonlinear term, we apply Lemma \ref{lemma:Gamma-R3} together with the $L^3-L^6-L^2$ estimate to obtain
        \begin{align}\label{Gamma_Pp f}
            ( \Gamma(f,f), \mathbf{P}^{\perp}f )
            \lesssim& \sum_{m\leq 2} 
            \| w_{\ell-m}\nabla_v^m f \|_{L^3_xL^2_v}
            \| f\|_{L^6_xL_v^2} \|\Pp f\| \notag\\
            \lesssim& \sum_{m\leq2,n\leq1} \|w_{\ell-m} \nabla^n_x\nabla^m_v f \| \left\| \nabla_x f\right\| \|\Pp f\| \notag\\
            \lesssim& \epsilon_0^2 \left\|\nabla_xf\right\|^2+\eta  \|\Pp f \|^2.
        \end{align}
        Substituting these estimates into \eqref{Pp f energy}, together with Young's inequality 
        \begin{align}
            \|\mathbf{P}^{\perp}f\|^2
            \lesssim
            \|\mathbf{P}^{\perp}f\langle v\rangle^{\frac{\gamma}{2}}\|^2
            +\eta \|\mathbf{P}^{\perp}f\langle v\rangle^{-\frac{\gamma}{2}}\|^2, \notag
        \end{align}
         we obtain
        \begin{align}
            &\frac{\mathrm{d}}{\mathrm{d}t} \|\mathbf{P}^{\perp}f \|^2+ \lambda \|\mathbf{P}^{\perp}f \|^2 
            \lesssim  \eta \|\mathbf{P}^{\perp}f\langle v\rangle^{-\frac{\gamma}{2}} \|^2 
            +\left\|\nabla_xf\right\|^2.\notag
        \end{align}
        Integrating it with respect to $t$ yields that
        \begin{align}
            \|\mathbf{P}^{\perp}f \|^2
            \lesssim& 
            e^{-\lambda t} \|\mathbf{P}^{\perp}f_0 \|^2
            +\eta\int_0^t e^{-\lambda(t-\tau)} \|\mathbf{P}^{\perp}f \langle v\rangle^{-\frac{\gamma}{2}} \|^2\mathrm{d}\tau
            +\int_0^t e^{-\lambda(t-\tau)} \|\nabla_x f \|^2\mathrm{d}\tau \notag\\
            \lesssim& 
            (1+t)^{-\frac52} \|\mathbf{P}^{\perp}f_0 \|^2
            +\eta \sup_t\left\{ (1+t)^\frac52\|\mathbf{P}^{\perp}f \langle v\rangle^{-\frac{\gamma}{2}} \|^2\right\}
            +\sup_t\left\{ (1+t)^\frac52\|\nabla_xf \|^2\right\}, \notag
        \end{align}
        where we have used the fact that
        \begin{align}
            \int_0^t e^{-\lambda (t-\tau)}(1+\tau)^{-\frac{5}{2}}\rmd \tau 
            \lesssim& e^{-\frac{\lambda}{2} t}\int_0^\frac{t}{2} (1+\tau)^{-\frac{5}{2}}\rmd \tau + \left(1+\frac{t}2\right)^{-\frac52} \int_\frac{t}{2}^t e^{-\lambda (t-\tau)}\rmd \tau 
            \lesssim (1+t)^{-\frac52}. \notag
        \end{align}
        This completes the proof of Lemma \ref{lemma Pp f}.
    \end{proof}
  \subsection{Microscopic part with weight}
    To control the weighted microscopic norm on the right-hand side of Lemma \ref{lemma Pp f}, we establish a hierarchy of time-weighted estimates for $\mathbf{P}^{\perp}f$ with increasing velocity weights.
  
    \begin{lemma}\label{lemma w Pp f}
        Let $-3<\gamma<0$, $N\geq4$ and ${\ell_0} \geq N$. Let $f(t,x,v)$ be the solution to the Cauchy problem \eqref{f-perturbation} and \eqref{f0}. Then it holds that
        \begin{align}
            &\sup_t\left\{ (1+t)^\frac52  \|w_{{\ell_0} }\Pp f\|^2 \right\}
            +\sup_t\left\{ (1+t)^\frac32  \|w_{{\ell_0}+\frac12 }\Pp f\|^2 \right\} \notag\\
            &+\sup_t\left\{ (1+t)^\frac12  \|w_{{\ell_0}+1 }\Pp f\|^2 \right\} 
            +\sup_t \|w_{{\ell_0}+\frac32 }\Pp f \|^2 \notag\\
            \lesssim & \|w_{{\ell_0}+\frac 32} f_0\|^2        
            +\sup_t\left\{(1+t)^\frac52\|\nabla_x f\|^2\right\}
            +\sup_t \left\{ (1+t)^\frac52  \|\Pp f \|^2 \right\}\notag\\
            &
            +(\eta+\epsilon_0^2)  \sum_{j\in \mathbb{Z}^+}\sup_t \left\{ (1+t)^\frac52  \| w_{{\ell_0} }\Pp f \|^2_{L^2_v(\mathbb{V}_j)L^2_x(\mathbb{R}^3)} \right\} \notag\\
            &+(\eta+\epsilon_0^2)  \sum_{j\in \mathbb{Z}^+}\sup_t \left\{ (1+t)^\frac32  \| w_{{\ell_0}+\frac12 }\Pp f \|^2_{L^2_v(\mathbb{V}_j)L^2_x(\mathbb{R}^3)} \right\} \notag\\
            &+(\eta+\epsilon_0^2)  \sum_{j\in \mathbb{Z}^+}\sup_t \left\{ (1+t)^\frac12  \| w_{{\ell_0}+1 }\Pp f \|^2_{L^2_v(\mathbb{V}_j)L^2_x(\mathbb{R}^3)} \right\} \notag\\
            &
            +(\eta+\epsilon_0^2)  \sum_{j\in \mathbb{Z}^+}\sup_t \| w_{{\ell_0}+\frac32}\Pp f \|^2_{L^2_v(\mathbb{V}_j)L^2_x(\mathbb{R}^3)} ,\notag
        \end{align}
        where $\eta$ and $\epsilon_0$ are chosen to be sufficiently small.
    \end{lemma}
    \begin{proof}
        Multiplying the microscopic equation \eqref{micro-eq} by $(1+t)^\frac52 w_{{\ell_0} }^2\Pp f$, we obtain
        \begin{align}
            &\frac{1}{2}\partial_t \{ (1+t)^\frac52 |w_{{\ell_0} }\Pp f|^2 \}
            +\frac{1}{2}\nu(v)(1+t)^\frac52 |w_{{\ell_0} }\Pp f|^2
            +\frac{1}{2}\nu(v)(1+t)^\frac52 |w_{{\ell_0} }\Pp f|^2 \notag\\
            =&-(1+t)^\frac52 (v\cdot\nabla_x\Pp f) w_{{\ell_0} }^2 \Pp f
            -(1+t)^\frac52 (v\cdot\nabla_x\P f) w_{{\ell_0} }^2 \Pp f
            +(1+t)^\frac52 \P(v\cdot\nabla_x f) w_{{\ell_0} }^2 \Pp f \notag\\
            &+(1+t)^\frac52 K\Pp f  w_{{\ell_0} }^2 \Pp f
            +(1+t)^\frac52 \Gamma(f,f) w_{{\ell_0} }^2 \Pp f
            +\frac54 (1+t)^\frac32 |w_{{\ell_0} } \Pp f|^2.\notag
        \end{align}
        Solving the above equation in time yields that
        \begin{align}
            &(1+t)^\frac52 |w_{{\ell_0} }\Pp f|^2
            +\int_0^t  e^{-\nu(v)(t-\tau)}\nu(v)(1+\tau)^\frac52|w_{{\ell_0} }\Pp f|^2\rmd \tau \notag\\
            \lesssim & e^{-\nu(v)t} |w_{{\ell_0} }\Pp f_0|^2
            -\int_0^t e^{-\nu(v)(t-\tau)} (1+\tau)^\frac52 (v\cdot\nabla_x\Pp f) w_{{\ell_0} }^2 \Pp f\rmd\tau \notag\\
            &-\int_0^t e^{-\nu(v)(t-\tau)} (1+\tau)^\frac52 (v\cdot\nabla_x\P f) w_{{\ell_0} }^2 \Pp f\rmd \tau
            +\int_0^t e^{-\nu(v)(t-\tau)}(1+\tau)^\frac52 \P(v\cdot\nabla_x f) w_{{\ell_0} }^2 \Pp f\rmd \tau \notag\\
            &+\int_0^t e^{-\nu(v)(t-\tau)}(1+\tau)^\frac52 K\Pp f  w_{{\ell_0} }^2 \Pp f\rmd\tau
            +\int_0^t e^{-\nu(v)(t-\tau)}(1+\tau)^\frac52  \Gamma(f,f) w_{{\ell_0} }^2 \Pp f\rmd\tau \notag\\
            &+\int_0^t e^{-\nu(v)(t-\tau)} (1+\tau)^\frac32 |w_{{\ell_0} } \Pp f|^2 \rmd\tau.\label{solve w Ppf}
        \end{align}
        Integrating \eqref{solve w Ppf} over $\mathbb{R}_v^3\times\mathbb{R}_x^3$, we obtain
        \begin{align}
            &(1+t)^\frac52 \|w_{{\ell_0} }\Pp f\|^2
            +\int_0^t \int_{\mathbb{R}_v^3\times\mathbb{R}_x^3} e^{-\nu(v)(t-\tau)}\nu(v)(1+\tau)^\frac52|w_{{\ell_0} }\Pp f|^2\rmd v \rmd x\rmd \tau \notag\\
            \lesssim & \|w_{{\ell_0} }\Pp f_0\|^2
            -\underbrace{\int_0^t \int_{\mathbb{R}_v^3\times\mathbb{R}_x^3} e^{-\nu(v)(t-\tau)} (1+\tau)^\frac52 (v\cdot\nabla_x\P f) w_{{\ell_0} }^2 \Pp f\rmd v \rmd x\rmd \tau }_{J_1} \notag \\
            &+\underbrace{\int_0^t \int_{\mathbb{R}_v^3\times\mathbb{R}_x^3} e^{-\nu(v)(t-\tau)}(1+\tau)^\frac52 \P(v\cdot\nabla_x f) w_{{\ell_0} }^2 \Pp f\rmd v \rmd x\rmd \tau}_{J_2} \notag\\
            &+\underbrace{\int_0^t \int_{\mathbb{R}_v^3\times\mathbb{R}_x^3} e^{-\nu(v)(t-\tau)}(1+\tau)^\frac52 K\Pp f  w_{{\ell_0} }^2 \Pp f\rmd v \rmd x\rmd\tau}_{J_3} \notag\\
            &+\underbrace{\int_0^t \int_{\mathbb{R}_v^3\times\mathbb{R}_x^3}  e^{-\nu(v)(t-\tau)}(1+\tau)^\frac52 \Gamma(f,f) w_{{\ell_0} }^2 \Pp f\rmd v \rmd x\rmd\tau }_{J_4}\notag\\
            &+\underbrace{\int_0^t \int_{\mathbb{R}_v^3\times\mathbb{R}_x^3} e^{-\nu(v)(t-\tau)} (1+\tau)^\frac32 |w_{{\ell_0} } \Pp f|^2 \rmd v \rmd x\rmd\tau}_{J_5}, \notag
        \end{align}
        where we have used the fact
        \begin{align*}
            \int_{\mathbb R^3_x}
            (v\cdot\nabla_x\Pp f)w_{\ell_0}^2 \Pp f\rmd x
            =
            0.
        \end{align*}
        For $J_1$, we use Young's inequality to obtain
        \begin{align}
            J_1 
            \lesssim& \int_0^t \int_{\mathbb{R}_v^3\times\mathbb{R}_x^3} e^{-\nu(v)(t-\tau)}(1+\tau)^\frac52 |\nabla_x(a,b,c)| \mu^\delta w_{{\ell_0} }^2 \Pp f  \rmd v \rmd x\rmd \tau \notag\\
            \lesssim&\int_0^t \int_{\mathbb{R}_v^3\times\mathbb{R}_x^3}e^{-\nu(v)(t-\tau)}\nu(v)(1+\tau)^\frac52 \left|\nabla_x(a,b,c)\right|^2 \mu^{\delta'} \rmd v \rmd x\rmd \tau \notag\\
            & +\eta\int_0^t \int_{\mathbb{R}_v^3\times\mathbb{R}_x^3}e^{-\nu(v)(t-\tau)}\nu(v)(1+\tau)^\frac52  |w_{{\ell_0} } \Pp f|^2\rmd v \rmd x\rmd \tau \notag\\
            \lesssim& \int_{\mathbb{R}_v^3}\int_0^te^{-\nu(v)(t-\tau)}\nu(v)\rmd \tau\mu^{\delta'} \rmd v \sup_\tau \left\{ (1+\tau)^\frac52  \left|\nabla_x(a,b,c)\right|^2_{L^2_x} \right\} \notag\\
            & +\eta\int_0^t \int_{\mathbb{R}_v^3\times\mathbb{R}_x^3}e^{-\nu(v)(t-\tau)}\nu(v)(1+\tau)^\frac52   |w_{{\ell_0} } \Pp f|^2\rmd v \rmd x\rmd \tau \notag\\
            \lesssim& \sup_t\left\{(1+t)^\frac52\|\nabla_x f\|^2\right\}
            +\eta\int_0^t \int_{\mathbb{R}_v^3\times\mathbb{R}_x^3}e^{-\nu(v)(t-\tau)}\nu(v)(1+\tau)^\frac52   |w_{{\ell_0} } \Pp f|^2\rmd v \rmd x\rmd \tau. \notag
        \end{align}
        Here we have used the fact that the polynomial factor in $v$ can be absorbed into the Maxwellian. 
        To estimate $J_2$, we recall the basis functions $\phi_i$ defined in \eqref{0phi-i}. By the definition of $\mathbf{P}$ and the Cauchy--Schwarz inequality in $v$, we have
        \begin{align}
            |\P(v\cdot\nabla_xf)|
            \lesssim
            \int_{\mathbb{R}_v^3} |v| |\nabla_x f||\phi_i|\rmd v \mu^\delta  
            \lesssim |\nabla_x f|_{L^2_v} \mu^\delta.\notag
        \end{align}
        Estimating $J_2$ in the same way as $J_1$, we obtain
        \begin{align}
            J_2
            \lesssim& \int_0^t \int_{\mathbb{R}_v^3\times\mathbb{R}_x^3} e^{-\nu(v)(t-\tau)}(1+\tau)^\frac52 |\nabla_x f|_{L^2_v} \mu^\delta w_{{\ell_0} }^2 \Pp f  \rmd v \rmd x\rmd \tau \notag\\
            \lesssim& \int_0^t \int_{\mathbb{R}_v^3\times\mathbb{R}_x^3}e^{-\nu(v)(t-\tau)}\nu(v)(1+\tau)^\frac52  \left|\nabla_xf\right|^2_{L^2_v}\mu^{\delta'}\rmd v \rmd x\rmd \tau \notag\\
            &+\eta\int_0^t \int_{\mathbb{R}_v^3\times\mathbb{R}_x^3}e^{-\nu(v)(t-\tau)}\nu(v)(1+\tau)^\frac52   |w_{{\ell_0} } \Pp f|^2\rmd v \rmd x\rmd \tau \notag\\
            \lesssim& \sup_t\left\{(1+t)^\frac52\|\nabla_x f\|^2\right\}
            +\eta\int_0^t \int_{\mathbb{R}_v^3\times\mathbb{R}_x^3}e^{-\nu(v)(t-\tau)}\nu(v)(1+\tau)^\frac52   |w_{{\ell_0} } \Pp f|^2\rmd v \rmd x\rmd \tau.\notag
        \end{align}
        To estimate $J_3$, we use the dyadic decomposition of the velocity space introduced in \eqref{dyadic}
        \begin{align}
            J_3 
            =& \int_0^t\int_{\mathbb{V}_0\times\mathbb{R}_x^3}e^{-\nu(v)(t-\tau)}(1+\tau)^\frac52 w_{{\ell_0} }^2 K\Pp f   \Pp f\rmd v \rmd x\rmd\tau \notag\\
            &+\sum_{j\in \mathbb{Z}^+}\int_0^t \int_{\mathbb{V}_j\times\mathbb{R}_x^3} e^{-\nu(v)(t-\tau)}(1+\tau)^\frac52 w_{{\ell_0} }^2 K\Pp f  \Pp f\rmd v \rmd x\rmd\tau \notag\\
            =&J_{3,1}+J_{3,2}.\notag
        \end{align}
        On the bounded-velocity region $\mathbb{V}_0$, the collision frequency $\nu(v)$ is bounded. Lemma \ref{lemma:K-R3} and Young's inequality yield
        \begin{align}
            J_{3,1} 
            \lesssim & \int_0^t\int_{\mathbb{V}_0\times\mathbb{R}_x^3}e^{-c_\gamma(t-\tau)}(1+\tau)^\frac52 w_{{\ell_0} }^2 K\Pp f   \Pp f\rmd v \rmd x\rmd\tau \notag\\
            \lesssim & \sup_t \left\{ (1+t)^\frac52\int_{\mathbb{R}_v^3\times\mathbb{R}_x^3} w_{{\ell_0} }^2 K\Pp f   \Pp f\rmd v \rmd x \right\} \notag\\
            \lesssim & \sup_t \left\{ (1+t)^\frac52 \int_{\mathbb{R}^3_x} \left(\eta|w_{{\ell_0} }\Pp f|_{L^2_\nu}+|w_{{\ell_0} }\bar{\chi}_C(v)\Pp f|_{L^2_v}\right) |w_{{\ell_0} }\Pp f|_{L^2_\nu} \rmd x \right\} \notag\\
            \lesssim & \eta\sup_t \left\{ (1+t)^\frac52 \|w_{{\ell_0} }\Pp f\|^2 \right\}+\sup_t \left\{ (1+t)^\frac52 \|\Pp f\|^2 \right\}.\notag
        \end{align}
        On each dyadic annulus $\mathbb{V}_j$, since $2^j \leq |v| < 2^{j+1}$ in $\mathbb{V}_j$, we have $(1+2^{j+1})^\gamma\lesssim \nu(v) \lesssim (1+2^j)^\gamma$ and 
        \begin{align}
            \int_0^t \mathbf{1}_{v \in \mathbb{V}_j}e^{-\nu(v)(t-\tau)}\nu(v)\rmd\tau 
            \lesssim \int_0^t e^{-C_\nu(1+2^{j+1})^\gamma(t-\tau)} (1+2^{j})^\gamma <\infty. \label{int Vj}
        \end{align}
        Combining \eqref{int Vj} with the localized estimate in Lemma \ref{lemma:K-Vj}, we obtain
        \begin{align}
            J_{3,2}
            \lesssim & \sum_{j\in \mathbb{Z}^+}\int_0^t \int_{\mathbb{V}_j\times\mathbb{R}_x^3} e^{-\nu(v)(t-\tau)}\nu(v) \nu^{-1}(v)(1+\tau)^\frac52 w_{{\ell_0} }^2 K\Pp f  \Pp f\rmd v \rmd x\rmd\tau \notag\\
            \lesssim & \sum_{j\in \mathbb{Z}^+} \sup_t \left\{ (1+t)^\frac52\int_{\mathbb{V}_j\times\mathbb{R}_x^3}  \nu^{-1}(v) w_{{\ell_0} }^2 K\Pp f  \Pp f\rmd v \rmd x \right\} \notag\\
            \lesssim & 
            \sum_{j\in \mathbb{Z}^+}\frac{\eta_1}{2^{2j}} \sup_t\left\{ (1+t)^\frac52\|w_{{\ell_0} }\Pp f\|^2 \right\}
			+\sum_{j\in \mathbb{Z}^+}\frac{1}{2^{2j}} \sup_t\left\{ (1+t)^\frac52\|\Pp f\|^2 \right\} \notag\\
			&+\eta_2\sum_{j\in \mathbb{Z}^+} \sup_t\left\{ (1+t)^\frac52 \| w_{{\ell_0} }\Pp f\|^2_{L^2_v(\widetilde{\mathbb{V}}_j)L_x^2(\mathbb{R}^3)}\right\}\notag\\
			&+\eta_3 \sum_{j\in \mathbb{Z}^+}\sup_t\left\{ (1+t)^\frac52 \| w_{{\ell_0} }\Pp f\|^2_{L^2_v({\mathbb{V}}_j)L_x^2(\mathbb{R}^3)}\right\} \notag\\
            \lesssim & (\eta_1+\eta_2) \sup_t \left\{ (1+t)^\frac52 \| w_{{\ell_0} }\Pp f\|^2 \right\}
            +\sup_t \left\{ (1+t)^\frac52 \| \Pp f\|^2 \right\} \notag\\
            &+(\eta_2+\eta_3)\sum_{j\in \mathbb{Z}^+}\sup_t \left\{ (1+t)^\frac52 \| w_{{\ell_0} }\Pp f\|^2_{L^2_v(\mathbb{V}_j)L^2_x(\mathbb{R}^3)} \right\}.\notag
        \end{align}
        We treat the nonlinear term $J_4$ by the same velocity decomposition
        \begin{align}
            J_4 
            =& \int_0^t\int_{\mathbb{V}_0\times\mathbb{R}_x^3}e^{-\nu(v)(t-\tau)}(1+\tau)^\frac52 w_{{\ell_0} }^2 \Gamma(f,f)  \Pp f\rmd v \rmd x\rmd\tau \notag\\
            &+\sum_{j\in \mathbb{Z}^+}\int_0^t \int_{\mathbb{V}_j\times\mathbb{R}_x^3} e^{-\nu(v)(t-\tau)}(1+\tau)^\frac52 w_{{\ell_0} }^2 \Gamma(f,f)  \Pp f\rmd v \rmd x\rmd\tau \notag\\
            =&J_{4,1}+J_{4,2}.\notag
        \end{align}
        Since $v$ is bounded on $\mathbb V_0$, the estimate \eqref{Gamma_Pp f} gives
        \begin{align}
            J_{4,1} 
            \lesssim & \int_0^t\int_{\mathbb{V}_0\times\mathbb{R}_x^3}e^{-C_\gamma(t-\tau)}(1+\tau)^\frac52   w_{{\ell_0}}^2\Gamma(f,f)  \Pp f\rmd v \rmd x\rmd\tau \notag\\
            \lesssim & \sup_t \left\{ (1+t)^\frac52\int_{\mathbb{R}_v^3\times\mathbb{R}_x^3} \Gamma(f,f)   \Pp f\rmd v \rmd x \right\} \notag\\
            \lesssim & \epsilon_0^2\sup_t \left\{ (1+t)^\frac52 \|\nabla_xf\|^2 \right\}
            +\eta \sup_t \left\{ (1+t)^\frac52 \|{\bf P}^\perp f \|^2 \right\}.\notag
        \end{align}
        On each dyadic annulus $\mathbb{V}_j$, we apply Lemma \ref{lemma:Gamma-Vj} and Young's inequality to obtain
        \begin{align}
            &\int_{\mathbb{V}_j\times\mathbb{R}_x^3}  \nu^{-1}(v) w_{{\ell_0} }^2 \Gamma(f,f)  \Pp f\rmd v \rmd x \notag\\
            \lesssim& \int_{\mathbb{R}_x^3}   \sum_{m\leq 2} |w_{{\ell_0}-m}  \nabla^m_v f|_{L_v^2(\mathbb{R}^3)} \left( |w_{{\ell_0} }{\bf P}f|_{L_v^2(\widehat{\mathbb{V}}_j)}+|w_{{\ell_0} }{\bf P}^\perp f|_{L_v^2(\widehat{\mathbb{V}}_j)} \right) |w_{{\ell_0} }{\bf P}^\perp f|_{L_v^2(\mathbb{V}_j)} \rmd x \notag\\
            \lesssim&  \sum_{m\leq 2} \| w_{{\ell_0}-m} \nabla^m_v f  \|_{L^3_x(\mathbb{R}^3)L_v^2(\mathbb{R}^3)} \notag\\
            &\times\left(  \mathbf{1}_{ \{1\leq j \leq 2 \}} \|w_{{\ell_0} }{\bf P}f \|_{L^6_x(\mathbb{R}^3)L_v^2(\widehat{\mathbb{V}}_j)}+ \mathbf{1}_{ \{j \geq 3\} }2^{-j}\|vw_{{\ell_0} }{\bf P}f \|_{L^6_x(\mathbb{R}^3)L_v^2(\widehat{\mathbb{V}}_j)} \right) \|w_{{\ell_0} }{\bf P}^\perp f \|_{L^2_x(\mathbb{R}^3)L_v^2(\mathbb{V}_j)} \notag\\
            &+ \sum_{m \leq 2} \| w_{{\ell_0}-m} \nabla^m_v f  \|_{L^\infty_x(\mathbb{R}^3)L_v^2(\mathbb{R}^3)} \|w_{{\ell_0} }{\bf P}^\perp f\|_{L^2_x(\mathbb{R}^3)L_v^2(\widehat{\mathbb{V}}_j)} \|w_{{\ell_0} }{\bf P}^\perp f \|_{L^2_x(\mathbb{R}^3)L_v^2(\mathbb{V}_j)}\notag\\
            \lesssim& \sum_{m\leq 2,n\leq 1} \| w_{{\ell_0}-m} \nabla^n_x \nabla^m_v f  \| 
            \frac{1}{2^{j}}|\nabla_x(a,b,c)|_{L^2_x(\mathbb{R}^3)}
            \|w_{{\ell_0} }{\bf P}^\perp f \|_{L^2_x(\mathbb{R}^3)L_v^2(\mathbb{V}_j)} \notag\\
            &+ \sum_{m,n\leq 2} \| w_{{\ell_0}-m} \nabla^n_x \nabla^m_v f  \|  \|w_{{\ell_0} }{\bf P}^\perp f \|_{L^2_x(\mathbb{R}^3)L_v^2(\widehat{\mathbb{V}}_j)} \|w_{{\ell_0} }{\bf P}^\perp f \|_{L^2_x(\mathbb{R}^3)L_v^2(\mathbb{V}_j)}\notag\\
            \lesssim&  \frac{\epsilon_0^2}{2^{2j}}  \|\nabla_xf \|^2 
            +\epsilon_0^2  \|w_{{\ell_0} }{\bf P}^\perp f \|_{L_v^2(\widehat{\mathbb{V}}_j)L^2_x(\mathbb{R}^3)}^2
            +\eta  \|w_{{\ell_0} }{\bf P}^\perp f \|^2_{L_v^2(\mathbb{V}_j)L^2_x(\mathbb{R}^3)},
            \label{Vj R3 w Gamma}
        \end{align}
        which implies
        \begin{align}
            J_{4,2}
            \lesssim & \sum_{j\in \mathbb{Z}^+}\int_0^t \int_{\mathbb{V}_j\times\mathbb{R}_x^3} e^{-\nu(v)(t-\tau)}\nu(v) \nu^{-1}(v)(1+\tau)^\frac52 w_{{\ell_0} }^2 \Gamma(f,f)  \Pp f\rmd v \rmd x\rmd\tau \notag\\
            \lesssim & \sum_{j\in \mathbb{Z}^+} \sup_t \left\{ (1+t)^\frac52\int_{\mathbb{V}_j\times\mathbb{R}_x^3}  \nu^{-1}(v) w_{{\ell_0} }^2 \Gamma(f,f)  \Pp f\rmd v \rmd x \right\} \notag\\
            \lesssim & \epsilon_0^2 \sup_t \left\{ (1+t)^\frac52  \|\nabla_xf\|^2 \right\}
            +\epsilon_0^2 \sup_t \left\{ (1+t)^\frac52  \| w_{{\ell_0} }\Pp f \|^2 \right\}\notag\\
            &+(\eta+\epsilon_0^2)  \sum_{j\in \mathbb{Z}^+}\sup_t \left\{ (1+t)^\frac52  \| w_{{\ell_0} }\Pp f \|^2_{L^2_v(\mathbb{V}_j)L^2_x(\mathbb{R}^3)} \right\}.
            \notag
        \end{align}
        Finally, using $\nu(v)\sim\langle v\rangle^\gamma$, we rewrite $J_5$ as
        \begin{align}
            J_5 \lesssim \int_0^t \int_{\mathbb{R}_v^3\times\mathbb{R}_x^3} e^{-\nu(v)(t-\tau)}\nu(v) (1+\tau)^\frac32 |w_{{\ell_0}+\frac12 } \Pp f|^2 \rmd v \rmd x\rmd\tau.\notag
        \end{align}
        Collecting the estimates for $J_1,\cdots,J_5$ and choosing $\eta$ and $\epsilon_0$ sufficiently small, we obtain
        \begin{align}
            &(1+t)^\frac52  \|w_{{\ell_0} }\Pp f \|^2
            +\int_0^t \int_{\mathbb{R}_v^3\times\mathbb{R}_x^3} e^{-\nu(v)(t-\tau)}\nu(v)(1+\tau)^\frac52|w_{{\ell_0} }\Pp f|^2\rmd v \rmd x\rmd \tau \notag\\
            \lesssim &  \|w_{{\ell_0} }\Pp f_0 \|^2           
            +\sup_t\left\{(1+t)^\frac52\|\nabla_x f\|^2\right\}
            +\sup_t \left\{ (1+t)^\frac52  \|\Pp f\|^2  \right\}
            \notag\\
            &+(\eta+\epsilon_0^2)  \sum_{j\in \mathbb{Z}^+}\sup_t \left\{ (1+t)^\frac52  \| w_{{\ell_0} }\Pp f \|^2_{L^2_v(\mathbb{V}_j)L^2_x(\mathbb{R}^3)} \right\}\notag\\
            &+\int_0^t \int_{\mathbb{R}_v^3\times\mathbb{R}_x^3} e^{-\nu(v)(t-\tau)}\nu(v) (1+\tau)^\frac32 |w_{{\ell_0} +\frac12} \Pp f|^2 \rmd v \rmd x\rmd\tau.\label{decay 5/2}
        \end{align}

        We next derive the estimate at the $(1+t)^{3/2}$ level by using $(1+t)^{3/2}w_{\ell_0+1/2}^2\mathbf{P}^{\perp}f$ as the multiplier. Repeating the preceding argument, we obtain
        \begin{align}
            &(1+t)^\frac32  \|w_{{\ell_0}+\frac12 }\Pp f \|^2
            +\int_0^t \int_{\mathbb{R}_v^3\times\mathbb{R}_x^3} e^{-\nu(v)(t-\tau)}\nu(v)(1+\tau)^\frac32|w_{{\ell_0}+\frac12 }\Pp f|^2\rmd v \rmd x\rmd \tau \notag\\
            \lesssim &  \|w_{{\ell_0}+\frac12 }\Pp f_0 \|^2            
            +\sup_t\left\{(1+t)^\frac32\|\nabla_x f\|^2\right\}
            +\sup_t \left\{ (1+t)^\frac32  \|\Pp f \|^2 \right\}\notag\\
            &
            +(\eta+\epsilon_0^2) \sum_{j\in \mathbb{Z}^+}\sup_t \left\{ (1+t)^\frac32  \| w_{{\ell_0}+\frac12 }\Pp f \|^2_{L^2_v(\mathbb{V}_j)L^2_x(\mathbb{R}^3)} \right\}\notag\\
            &+\int_0^t \int_{\mathbb{R}_v^3\times\mathbb{R}_x^3} e^{-\nu(v)(t-\tau)}\nu(v) (1+\tau)^\frac12 |w_{{\ell_0}+1 } \Pp f|^2 \rmd v \rmd x\rmd\tau. \label{decay 3/2}
        \end{align}
        Proceeding similarly with the multiplier $(1+t)^{1/2}w_{\ell_0+1}^2\mathbf{P}^{\perp}f$, we obtain
        \begin{align}
            &(1+t)^\frac12  \|w_{{\ell_0}+1 }\Pp f \|^2
            +\int_0^t \int_{\mathbb{R}_v^3\times\mathbb{R}_x^3} e^{-\nu(v)(t-\tau)}\nu(v)(1+\tau)^\frac12|w_{{\ell_0}+1 }\Pp f|^2\rmd v \rmd x\rmd \tau \notag\\
            \lesssim &  \|w_{{\ell_0}+1 }\Pp f_0 \|^2            
            +\sup_t\left\{(1+t)^\frac12\|\nabla_x f\|^2\right\}
            +\sup_t \left\{ (1+t)^\frac12  \|\Pp f \|^2 \right\}\notag\\
            &
            +(\eta+\epsilon_0^2)  \sum_{j\in \mathbb{Z}^+}\sup_t \left\{ (1+t)^\frac12  \| w_{{\ell_0}+1 }\Pp f \|^2_{L^2_v(\mathbb{V}_j)L^2_x(\mathbb{R}^3)} \right\}\notag\\
            &+\int_0^t \int_{\mathbb{R}_v^3\times\mathbb{R}_x^3} e^{-\nu(v)(t-\tau)} \nu(v) (1+\tau)^{-\frac12} |w_{{\ell_0}+\frac 32 } \Pp f|^2 \rmd v \rmd x\rmd\tau. \label{decay 1/2}
        \end{align}
        To control the last term in \eqref{decay 1/2}, we multiply \eqref{micro-eq} by $w_{{\ell_0}+\frac32}^2\Pp f$ and argue as above to obtain
        \begin{align}
            & \|w_{{\ell_0}+\frac32 }\Pp f \|^2
            +\int_0^t \int_{\mathbb{R}_v^3\times\mathbb{R}_x^3} e^{-\nu(v)(t-\tau)}\nu(v)|w_{{\ell_0}+\frac32}\Pp f|^2\rmd v \rmd x\rmd \tau \notag\\
            \lesssim &  \|w_{{\ell_0}+\frac32 }\Pp f_0 \|^2   
            +\sup_t\|\nabla_x f\|^2
            +\sup_t\|\Pp f \|^2 \notag\\
            &+(\eta+\epsilon_0^2)  \sum_{j\in \mathbb{Z}^+}\sup_t \| w_{{\ell_0}+\frac32}\Pp f \|^2_{L^2_v(\mathbb{V}_j)L^2_x(\mathbb{R}^3)}. \label{decay -1/2}
        \end{align}
        Taking the supremum over $t$ in \eqref{decay -1/2} gives
        \begin{align}
            \sup_t \|w_{{\ell_0}+\frac32 }\Pp f \|^2 
            \lesssim &  \|w_{{\ell_0}+\frac32 }f_0 \|^2   
            +\sup_t\|\nabla_x f\|^2
            +\sup_t\|\Pp f \|^2 \notag\\
            &+(\eta+\epsilon_0^2)  \sum_{j\in \mathbb{Z}^+}\sup_t \| w_{{\ell_0}+\frac32}\Pp f \|^2_{L^2_v(\mathbb{V}_j)L^2_x(\mathbb{R}^3)}.\notag
        \end{align}
        Combining \eqref{decay 1/2} and \eqref{decay -1/2}, we obtain the estimate at the $(1+t)^{1/2}$ level
        \begin{align}
            \sup_t\left\{(1+t)^\frac12  \|w_{{\ell_0}+1 }\Pp f \|^2\right\}
            \lesssim & \|w_{{\ell_0}+\frac 32} f_0 \|^2       
            +\sup_t\left\{(1+t)^\frac12\|\nabla_x f\|^2\right\}
            +\sup_t \left\{ (1+t)^\frac12  \|\Pp f \|^2 \right\}\notag\\
            & 
            +(\eta+\epsilon_0^2)  \sum_{j\in \mathbb{Z}^+}\sup_t \left\{ (1+t)^\frac12  \| w_{{\ell_0}+1 }\Pp f \|^2_{L^2_v(\mathbb{V}_j)L^2_x(\mathbb{R}^3)} \right\} \notag\\
            &
            +(\eta+\epsilon_0^2)  \sum_{j\in \mathbb{Z}^+}\sup_t \| w_{{\ell_0}+\frac32}\Pp f \|^2_{L^2_v(\mathbb{V}_j)L^2_x(\mathbb{R}^3)} .\notag
        \end{align}
        Using \eqref{decay 3/2} together with \eqref{decay 1/2} and \eqref{decay -1/2}, we similarly obtain the estimate at the $(1+t)^{3/2}$ level
        \begin{align}
            \sup_t\left\{(1+t)^\frac32  \|w_{{\ell_0}+\frac12 }\Pp f \|^2\right\}
            \lesssim & \| w_{{\ell_0}+\frac 32} f_0 \|^2          
            +\sup_t\left\{(1+t)^\frac32\|\nabla_x f\|^2\right\}
            +\sup_t \left\{ (1+t)^\frac32  \|\Pp f \|^2 \right\}\notag\\
            &
            +(\eta+\epsilon_0^2)  \sum_{j\in \mathbb{Z}^+}\sup_t \left\{ (1+t)^\frac32  \| w_{{\ell_0}+\frac12 }\Pp f \|^2_{L^2_v(\mathbb{V}_j)L^2_x(\mathbb{R}^3)} \right\}\notag\\
            &
            +(\eta+\epsilon_0^2)  \sum_{j\in \mathbb{Z}^+}\sup_t \left\{ (1+t)^\frac12  \| w_{{\ell_0}+1 }\Pp f \|^2_{L^2_v(\mathbb{V}_j)L^2_x(\mathbb{R}^3)} \right\}\notag\\
            &
            +(\eta+\epsilon_0^2)  \sum_{j\in \mathbb{Z}^+}\sup_t \| w_{{\ell_0}+\frac32}\Pp f \|^2_{L^2_v(\mathbb{V}_j)L^2_x(\mathbb{R}^3)} . \notag
        \end{align}
        Finally, combining \eqref{decay 5/2}, \eqref{decay 3/2}, \eqref{decay 1/2} and \eqref{decay -1/2}, we obtain the hierarchy at the $(1+t)^{5/2}$ level
        \begin{align}
            \sup_t\left\{ (1+t)^\frac52  \|w_{{\ell_0} }\Pp f \|^2\right\}
            \lesssim &  \|w_{{\ell_0}+\frac 32} f_0 \|^2         
            +\sup_t\left\{(1+t)^\frac52\|\nabla_x f\|^2\right\}
            +\sup_t \left\{ (1+t)^\frac52  \|\Pp f\|^2  \right\}\notag\\
            &
            +(\eta+\epsilon_0^2)  \sum_{j\in \mathbb{Z}^+}\sup_t \left\{ (1+t)^\frac52  \| w_{{\ell_0} }\Pp f \|^2_{L^2_v(\mathbb{V}_j)L^2_x(\mathbb{R}^3)} \right\} \notag\\
            &+(\eta+\epsilon_0^2)  \sum_{j\in \mathbb{Z}^+}\sup_t \left\{ (1+t)^\frac32  \| w_{{\ell_0}+\frac12 }\Pp f \|^2_{L^2_v(\mathbb{V}_j)L^2_x(\mathbb{R}^3)} \right\} \notag\\
            &+(\eta+\epsilon_0^2)  \sum_{j\in \mathbb{Z}^+}\sup_t \left\{ (1+t)^\frac12  \| w_{{\ell_0}+1 }\Pp f \|^2_{L^2_v(\mathbb{V}_j)L^2_x(\mathbb{R}^3)} \right\}\notag\\
            &
            +(\eta+\epsilon_0^2)  \sum_{j\in \mathbb{Z}^+}\sup_t \| w_{{\ell_0}+\frac32}\Pp f \|^2_{L^2_v(\mathbb{V}_j)L^2_x(\mathbb{R}^3)} .\notag
        \end{align}
        This completes the proof of Lemma \ref{lemma w Pp f}.
    \end{proof}
  \subsection{Microscopic part on velocity annulus $\mathbb{V}_j$}
    It remains to control the dyadic remainder terms appearing in Lemma \ref{lemma w Pp f}. We establish the required localized estimates on each velocity annulus $\mathbb{V}_j$.
    \begin{lemma}\label{lemma Vj}
    Let $-3<\gamma<0$, $N\geq4$ and ${\ell_0} \geq N$. Let $f(t,x,v)$ be the solution to the Cauchy problem \eqref{f-perturbation} and \eqref{f0}. Then it holds that
        \begin{align}
            &\sum_{j\in\mathbb{Z}^+} \sup_t \left\{ (1+t)^\frac52  \|w_{{\ell_0} }\Pp f \|^2_{L^2_v(\mathbb{V}_j)L^2_x(\mathbb{R}^3)} \right\}
            +\sum_{j\in\mathbb{Z}^+} \sup_t\left\{ (1+t)^\frac32  \|w_{{\ell_0}+\frac12 }\Pp f \|^2_{L^2_v(\mathbb{V}_j)L^2_x(\mathbb{R}^3)} \right\} \notag\\
            &+\sum_{j\in\mathbb{Z}^+} \sup_t \left\{ (1+t)^\frac12  \|w_{{\ell_0}+1 }\Pp f \|^2_{L^2_v(\mathbb{V}_j)L^2_x(\mathbb{R}^3)}\right\} 
            +\sum_{j\in \mathbb{Z}^+}\sup_t \| w_{{\ell_0}+\frac32}\Pp f \|^2_{L^2_v(\mathbb{V}_j)L^2_x(\mathbb{R}^3)} \notag\\
            \lesssim& 
            \| w_{{\ell_0}+\frac 32}f_0 \|^2
            +\sup_t\left\{(1+t)^\frac52 \|\nabla_x f\|^2\right\}
            +\sup_t \left\{ (1+t)^\frac52 \| \Pp f\|^2 \right\}\notag\\
            &+(\eta+\epsilon_0^2) \sup_t\left\{ (1+t)^\frac52\|w_{{\ell_0} }\Pp f\|^2 \right\} 
            +(\eta+\epsilon_0^2) \sup_t\left\{ (1+t)^\frac32\|w_{{\ell_0}+\frac12 }\Pp f\|^2 \right\} \notag\\
            &+(\eta+\epsilon_0^2) \sup_t\left\{ (1+t)^\frac12\|w_{{\ell_0}+1 }\Pp f\|^2 \right\} 
            +(\eta+\epsilon_0^2) \sup_t \|w_{{\ell_0}+\frac32 }\Pp f\|^2, \notag
        \end{align}
        where $\eta$ and $\epsilon_0$ are chosen to be sufficiently small.
    \end{lemma}
    \begin{proof}
        For each $j\geq1$, we integrate \eqref{solve w Ppf} over $\mathbb{V}_j\times\mathbb{R}_x^3$ to obtain the localized weighted energy inequality
        \begin{align}
            &(1+t)^\frac52  \|w_{{\ell_0} }\Pp f \|^2_{L^2_v(\mathbb{V}_j)L^2_x(\mathbb{R}^3)}
            +\int_0^t \int_{\mathbb{V}_j\times\mathbb{R}_x^3} e^{-\nu(v)(t-\tau)}\nu(v)(1+\tau)^\frac52|w_{{\ell_0} }\Pp f|^2\rmd v \rmd x\rmd \tau \notag\\
            \lesssim &  \|w_{{\ell_0} }\Pp f_0 \|^2_{L^2_v(\mathbb{V}_j)L^2_x(\mathbb{R}^3)}
            -\underbrace{\int_0^t \int_{\mathbb{V}_j\times\mathbb{R}_x^3} e^{-\nu(v)(t-\tau)} (1+\tau)^\frac52 (v\cdot\nabla_x\P f) w_{{\ell_0} }^2 \Pp f\rmd v \rmd x\rmd \tau }_{M_1} \notag \\
            &+\underbrace{\int_0^t \int_{\mathbb{V}_j\times\mathbb{R}_x^3} e^{-\nu(v)(t-\tau)}(1+\tau)^\frac52 \P(v\cdot\nabla_x f) w_{{\ell_0} }^2 \Pp f\rmd v \rmd x\rmd \tau}_{M_2} \notag\\
            &+\underbrace{\int_0^t \int_{\mathbb{V}_j\times\mathbb{R}_x^3} e^{-\nu(v)(t-\tau)}(1+\tau)^\frac52 K\Pp f  w_{{\ell_0} }^2 \Pp f\rmd v \rmd x\rmd\tau}_{M_3} \notag\\
            &+\underbrace{\int_0^t \int_{\mathbb{V}_j\times\mathbb{R}_x^3}  e^{-\nu(v)(t-\tau)}(1+\tau)^\frac52 \Gamma(f,f) w_{{\ell_0} }^2 \Pp f\rmd v \rmd x\rmd\tau }_{M_4}\notag\\
            &+\underbrace{\int_0^t \int_{\mathbb{V}_j\times\mathbb{R}_x^3} e^{-\nu(v)(t-\tau)} (1+\tau)^\frac32 |w_{{\ell_0} } \Pp f|^2 \rmd v \rmd x\rmd\tau}_{M_5}.\notag
        \end{align}
        We first estimate $M_1$. Using \eqref{int Vj} and Young's inequality, we obtain
        \begin{align}
            M_1
            \lesssim& \sup_t \left\{ (1+t)^\frac52 \int_{\mathbb{V}_j\times\mathbb{R}_x^3}   \nu^{-1}(v) (v\cdot\nabla_x\P f) w_{{\ell_0} }^2 \Pp f\rmd v \rmd x \right\} \notag\\
            \lesssim& \sup_t \left\{ (1+t)^\frac52 \int_{\mathbb{V}_j\times\mathbb{R}_x^3}   |\nabla_x(a,b,c)| |v|^{-1}\mu^\delta w_{{\ell_0} } \Pp f\rmd v \rmd x \right\} \notag\\
            \lesssim& \frac{1}{2^{2j}}\sup_t\left\{(1+t)^\frac52 \|\nabla_xf \|^2 \right\}
            +\eta\sup_t \left\{ (1+t)^\frac52  \|w_{{\ell_0} }\Pp f \|^2_{L^2_v(\mathbb{V}_j)L^2_x(\mathbb{R}^3)} \right\}.\notag
        \end{align}
        Here, as before, the polynomial factor is absorbed into the Maxwellian.
        The term $M_2$ is estimated in the same way.
        By \eqref{int Vj}, \eqref{0phi-i} and Young's inequality, we obtain
        \begin{align}
            M_2
            \lesssim& \sup_t \left\{ (1+t)^\frac52 \int_{\mathbb{V}_j\times\mathbb{R}_x^3}   \nu^{-1}(v) \P(v\cdot\nabla_x f) w_{{\ell_0} }^2 \Pp f\rmd v \rmd x \right\} \notag\\
            \lesssim& \sup_t \left\{ (1+t)^\frac52 \int_{\mathbb{V}_j\times\mathbb{R}_x^3}   |\nabla_x f|_{L^2_v} |v|^{-1} \mu^\delta  w_{{\ell_0} } \Pp f\rmd v \rmd x \right\} \notag\\
            \lesssim& \frac{1}{2^{2j}}\sup_t\left\{(1+t)^\frac52 \|\nabla_x f\|^2\right\}
            +\eta\sup_t \left\{ (1+t)^\frac52  \|w_{{\ell_0} }\Pp f \|^2_{L^2_v(\mathbb{V}_j)L^2_x(\mathbb{R}^3)} \right\}.\notag
        \end{align}
        For $M_3$, we combine \eqref{int Vj} with the localized estimate for $K$ in Lemma \ref{lemma:K-Vj} to obtain
        \begin{align}
            M_3
            \lesssim& \sup_t \left\{ (1+t)^\frac52 \int_{\mathbb{V}_j\times\mathbb{R}_x^3} \nu^{-1}(v) w_{{\ell_0} }^2 K\Pp f   \Pp f\rmd v \rmd x \right\} \notag\\
            \lesssim&
            \frac{\eta_1}{2^{2j}} \sup_t\left\{ (1+t)^\frac52\|w_{{\ell_0} }\Pp f\|^2 \right\}
            +\frac{1}{2^{2j}} \sup_t\left\{ (1+t)^\frac52\|\Pp f\|^2 \right\} \notag\\
            &+\eta_2\sup_t\left\{ (1+t)^\frac52 \| w_{{\ell_0} }\Pp f\|^2_{L^2(\widetilde{\mathbb{V}}_j)L_x^2(\mathbb{R}^3)}\right\}
            +\eta_3 \sup_t\left\{ (1+t)^\frac52 \| w_{{\ell_0} }\Pp f\|^2_{L^2_v({\mathbb{V}}_j)L_x^2(\mathbb{R}^3)}\right\}.\notag
        \end{align}
        For the nonlinear term $M_4$, the time-integral estimate \eqref{int Vj} and the localized bound \eqref{Vj R3 w Gamma} yield
        \begin{align}
            M_4
            \lesssim& \sup_t \left\{ (1+t)^\frac52 \int_{\mathbb{V}_j\times\mathbb{R}_x^3}   \nu^{-1}(v) w_{{\ell_0} }^2 \Gamma(f,f)   \Pp f\rmd v \rmd x \right\} \notag\\
            \lesssim& \frac{\epsilon_0^2}{2^{2j}}\sup_t \left\{ (1+t)^\frac52  \|\nabla_xf\|^2  \right\}
            +\epsilon_0^2 \sup_t \left\{ (1+t)^\frac52  \|w_{{\ell_0} }{\bf P}^\perp f \|_{L_v^2(\widehat{\mathbb{V}}_j)L^2_x(\mathbb{R}^3)}^2\right\}\notag\\
            &+\eta \sup_t \left\{ (1+t)^\frac52 \|w_{{\ell_0} }{\bf P}^\perp f \|^2_{L_v^2(\mathbb{V}_j)L^2_x(\mathbb{R}^3)}\right\}.\notag
        \end{align}
        Finally, using $\nu(v)\sim\langle v\rangle^\gamma$, we can express $M_5$ as
        \begin{align}
            M_5&\lesssim\int_0^t \int_{\mathbb{V}_j\times\mathbb{R}_x^3} e^{-\nu(v)(t-\tau)}\nu(v) (1+\tau)^\frac32 |w_{{\ell_0}+\frac12 } \Pp f|^2 \rmd v \rmd x\rmd\tau.\notag
        \end{align}
        Collecting the bounds for $M_1,\cdots,M_5$, we arrive at the localized weighted energy inequality
        \begin{align}
            &(1+t)^\frac52  \|w_{{\ell_0} }\Pp f \|^2_{L^2_v(\mathbb{V}_j)L^2_x(\mathbb{R}^3)}
            +\int_0^t \int_{\mathbb{V}_j\times\mathbb{R}_x^3} e^{-\nu(v)(t-\tau)}\nu(v)(1+\tau)^\frac52|w_{{\ell_0} }\Pp f|^2\rmd v \rmd x\rmd \tau  \notag\\
            \lesssim &   \|w_{{\ell_0} }\Pp f_0 \|^2_{L^2_v(\mathbb{V}_j)L^2_x(\mathbb{R}^3)}
            +\frac{1}{2^{2j}}\sup_t\left\{(1+t)^\frac52 \|\nabla_x f\|^2\right\}
            +\frac{1}{2^{2j}}\sup_t \left\{ (1+t)^\frac52  \| \Pp f\|^2  \right\}\notag\\
            &+\frac{\eta}{2^{2j}} \sup_t\left\{ (1+t)^\frac52\|w_{{\ell_0} }\Pp f\|^2 \right\}\notag\\
            &
            +\eta'\sup_t\left\{ (1+t)^\frac52 \| w_{{\ell_0} }\Pp f\|^2_{L^2(\widetilde{\mathbb{V}}_j)L_x^2(\mathbb{R}^3)}\right\}+\epsilon_0^2 \sup_t\left\{ (1+t)^\frac52 \| w_{{\ell_0} }\Pp f\|^2_{L^2(\widehat{\mathbb{V}}_j)L_x^2(\mathbb{R}^3)}\right\}\notag\\
            &+ \int_0^t \int_{\mathbb{V}_j\times\mathbb{R}_x^3} e^{-\nu(v)(t-\tau)}\nu(v) (1+\tau)^\frac32 |w_{{\ell_0}+\frac12 } \Pp f|^2 \rmd v \rmd x\rmd\tau .\label{Vj 5/2}
        \end{align}
        Repeating the preceding argument with the multiplier $(1+t)^{3/2}w_{\ell_0+1/2}^2\mathbf{P}^{\perp}f$, we solve the resulting differential inequality and integrate over $\mathbb{V}_j\times\mathbb{R}_x^3$ to obtain
        \begin{align}
            &(1+t)^\frac32  \|w_{{\ell_0}+\frac12 }\Pp f \|^2_{L^2_v(\mathbb{V}_j)L^2_x(\mathbb{R}^3)}
            +\int_0^t \int_{\mathbb{V}_j\times\mathbb{R}_x^3} e^{-\nu(v)(t-\tau)}\nu(v)(1+\tau)^\frac32|w_{{\ell_0}+\frac12 }\Pp f|^2\rmd v \rmd x\rmd \tau  \notag\\
            \lesssim &   \|w_{{\ell_0}+\frac12 }\Pp f_0 \|^2_{L^2_v(\mathbb{V}_j)L^2_x(\mathbb{R}^3)}
            +\frac{1}{2^{2j}}\sup_t\left\{(1+t)^\frac32 \|\nabla_x f\|^2\right\}
            +\frac{1}{2^{2j}}\sup_t \left\{ (1+t)^\frac32  \| \Pp f \|^2 \right\}\notag\\
            &
            +\frac{\eta}{2^{2j}} \sup_t\left\{ (1+t)^\frac32\|w_{{\ell_0}+\frac12 }\Pp f\|^2 \right\} \notag\\
            &+\eta'\sup_t\left\{ (1+t)^\frac32 \| w_{{\ell_0}+\frac12 }\Pp f\|^2_{L^2(\widetilde{\mathbb{V}}_j)L_x^2(\mathbb{R}^3)}\right\}
            +\epsilon_0^2 \sup_t\left\{ (1+t)^\frac32 \| w_{{\ell_0}+\frac12 }\Pp f\|^2_{L^2(\widehat{\mathbb{V}}_j)L_x^2(\mathbb{R}^3)}\right\}
             \notag\\
            &+\int_0^t \int_{\mathbb{V}_j\times\mathbb{R}_x^3} e^{-\nu(v)(t-\tau)}\nu(v) (1+\tau)^\frac12 |w_{{\ell_0}+1 } \Pp f|^2 \rmd v \rmd x\rmd\tau.\label{Vj 3/2}
        \end{align}
        The same procedure with the multiplier $(1+t)^\frac12 w_{{\ell_0}+1}^2\Pp f$ gives
        \begin{align}
            &(1+t)^\frac12  \|w_{{\ell_0}+1 }\Pp f \|^2_{L^2_v(\mathbb{V}_j)L^2_x(\mathbb{R}^3)}
            +\int_0^t \int_{\mathbb{V}_j\times\mathbb{R}_x^3} e^{-\nu(v)(t-\tau)}\nu(v)(1+\tau)^\frac12|w_{{\ell_0}+1 }\Pp f|^2\rmd v \rmd x\rmd \tau  \notag\\
            \lesssim &   \|w_{{\ell_0}+1 }\Pp f_0\|^2_{L^2_v(\mathbb{V}_j)L^2_x(\mathbb{R}^3)}
            +\frac{1}{2^{2j}}\sup_t\left\{(1+t)^\frac12 \|\nabla_x f\|^2\right\}
            +\frac{1}{2^{2j}}\sup_t \left\{ (1+t)^\frac12  \| \Pp f\|^2  \right\}
             \notag\\
            &
            +\frac{\eta}{2^{2j}} \sup_t\left\{ (1+t)^\frac12\|w_{{\ell_0}+1 }\Pp f\|^2 \right\} \notag\\
            &
            +\eta'\sup_t\left\{ (1+t)^\frac12 \| w_{{\ell_0}+1 }\Pp f\|^2_{L^2(\widetilde{\mathbb{V}}_j)L_x^2(\mathbb{R}^3)}\right\}
            +\epsilon_0^2\sup_t\left\{ (1+t)^\frac12 \| w_{{\ell_0}+1 }\Pp f\|^2_{L^2(\widehat{\mathbb{V}}_j)L_x^2(\mathbb{R}^3)}\right\}\notag\\
            &+\int_0^t \int_{\mathbb{V}_j\times\mathbb{R}_x^3} e^{-\nu(v)(t-\tau)}\nu(v) (1+\tau)^{-\frac12} |w_{{\ell_0}+\frac 32 } \Pp f|^2 \rmd v \rmd x\rmd\tau .\label{Vj 1/2}
        \end{align}
        To control the remaining term involving the weight $w_{\ell_0+3/2}$, we use $w_{\ell_0+3/2}^2\mathbf{P}^{\perp}f$ as the multiplier and obtain
        \begin{align}
            &\|w_{{\ell_0}+\frac32}\Pp f \|^2_{L^2_v(\mathbb{V}_j)L^2_x(\mathbb{R}^3)}
            +\int_0^t \int_{\mathbb{V}_j\times\mathbb{R}_x^3} e^{-\nu(v)(t-\tau)}\nu(v)|w_{{\ell_0}+\frac32 }\Pp f|^2\rmd v \rmd x\rmd \tau  \notag\\
            \lesssim &   \|w_{{\ell_0}+\frac32 }\Pp f_0\|^2_{L^2_v(\mathbb{V}_j)L^2_x(\mathbb{R}^3)}
            +\frac{1}{2^{2j}}\sup_t \|\nabla_x f\|^2
            +\frac{1}{2^{2j}}\sup_t \| \Pp f\|^2
            +\frac{\eta}{2^{2j}} \sup_t \|w_{{\ell_0}+\frac32 }\Pp f\|^2 \notag\\
            &
            +\eta'\sup_t \| w_{{\ell_0}+\frac32 }\Pp f\|^2_{L^2(\widetilde{\mathbb{V}}_j)L_x^2(\mathbb{R}^3)}
            +\epsilon_0^2\sup_t \| w_{{\ell_0}+\frac32 }\Pp f\|^2_{L^2(\widehat{\mathbb{V}}_j)L_x^2(\mathbb{R}^3)}.\label{Vj -1/2}
        \end{align}
        Taking the supremum over time in \eqref{Vj -1/2} and summing over $j\in\mathbb{Z}^+$, we obtain
        \begin{align}
            &\sum_{j\in\mathbb{Z}^+}\sup_t \|w_{{\ell_0}+\frac32}\Pp f \|^2_{L^2_v(\mathbb{V}_j)L^2_x(\mathbb{R}^3)} \notag\\
            \lesssim &  \|w_{{\ell_0}+\frac 32} f_0\|^2
            +\sup_t \|\nabla_x f\|^2 
            +\sup_t \| \Pp f \|^2 
            +(\eta+\epsilon_0^2) \sup_t\|w_{{\ell_0}+\frac32}\Pp f\|^2 .\notag
        \end{align}
        Combining \eqref{Vj 1/2} and \eqref{Vj -1/2}, we take the supremum over time and sum over $j\in\mathbb{Z}^+$ to obtain
        \begin{align}
            &\sum_{j\in\mathbb{Z}^+}\sup_t\left\{ (1+t)^\frac12  \|w_{{\ell_0}+1 }\Pp f \|^2_{L^2_v(\mathbb{V}_j)L^2_x(\mathbb{R}^3)} \right\} \notag\\
            \lesssim &  \|w_{{\ell_0}+\frac 32} f_0\|^2
            +\sup_t\left\{(1+t)^\frac12 \|\nabla_x f\|^2\right\}
            +\sup_t \left\{ (1+t)^\frac12  \| \Pp f \|^2 \right\}
            \notag\\
            &+(\eta+\epsilon_0^2) \sup_t\left\{ (1+t)^\frac12\|w_{{\ell_0}+1 }\Pp f\|^2 \right\}
            +(\eta+\epsilon_0^2) \sup_t \|w_{{\ell_0}+\frac32 }\Pp f\|^2.\notag
        \end{align}
        Similarly, combining \eqref{Vj 3/2}, \eqref{Vj 1/2} and \eqref{Vj -1/2}, we obtain
        \begin{align}
            &\sum_{j\in\mathbb{Z}^+}\sup_t\left\{ (1+t)^\frac32  \|w_{{\ell_0}+\frac12 }\Pp f \|^2_{L^2_v(\mathbb{V}_j)L^2_x(\mathbb{R}^3)} \right\} \notag\\
            \lesssim &  \|w_{{\ell_0}+\frac 32} f_0\|^2
            +\sup_t\left\{(1+t)^\frac32 \|\nabla_x f\|^2\right\}
            +\sup_t \left\{ (1+t)^\frac32  \| \Pp f \|^2 \right\}
            \notag\\
            &+(\eta+\epsilon_0^2) \sup_t\left\{ (1+t)^\frac32\|w_{{\ell_0}+\frac12 }\Pp f\|^2 \right\}
            +(\eta+\epsilon_0^2) \sup_t\left\{ (1+t)^\frac12\|w_{{\ell_0}+1 }\Pp f\|^2 \right\}
            \notag\\
            & +(\eta+\epsilon_0^2) \sup_t \|w_{{\ell_0}+\frac32 }\Pp f\|^2.\notag
        \end{align}
        Finally, combining \eqref{Vj 5/2}, \eqref{Vj 3/2}, \eqref{Vj 1/2} and \eqref{Vj -1/2}, we obtain
        \begin{align}
            &\sum_{j\in\mathbb{Z}^+}\sup_t\left\{ (1+t)^\frac52  \|w_{{\ell_0} }\Pp f \|^2_{L^2_v(\mathbb{V}_j)L^2_x(\mathbb{R}^3)} \right\}
            \notag\\
            \lesssim &  \|w_{{\ell_0}+\frac 32} f_0\|^2
            +\sup_t\left\{(1+t)^\frac52 \|\nabla_xf\|^2\right\}
            +\sup_t \left\{ (1+t)^\frac52 \| \Pp f\|^2 \right\}\notag\\
            &+(\eta+\epsilon_0^2) \sup_t\left\{ (1+t)^\frac52\|w_{{\ell_0} }\Pp f\|^2 \right\} 
            +(\eta+\epsilon_0^2) \sup_t\left\{ (1+t)^\frac32\|w_{{\ell_0}+\frac12 }\Pp f\|^2 \right\} 
            \notag\\
            &+(\eta+\epsilon_0^2) \sup_t\left\{ (1+t)^\frac12\|w_{{\ell_0}+1 }\Pp f\|^2 \right\}
            +(\eta+\epsilon_0^2) \sup_t \|w_{{\ell_0}+\frac32 }\Pp f\|^2.\notag
        \end{align}
        This completes the proof of Lemma \ref{lemma Vj}.
    \end{proof}

    \noindent{\bf Proof of Theorem \ref{thm}.}
    Combining Lemmas \ref{lemma w Pp f} and \ref{lemma Vj}, and choosing  $\eta$ and $\epsilon_0$ sufficiently small, we absorb the dyadic remainder terms and obtain the desired weighted estimate for $\mathbf{P}^{\perp}f$
    \begin{align}
        &\sup_t\left\{ (1+t)^\frac52  \|w_{{\ell_0} }\Pp f\|^2 \right\}
        \lesssim  \|w_{{\ell_0}+\frac 32} f_0\|^2          
        +\sup_t\left\{(1+t)^\frac52\|\nabla_x f\|^2\right\}
        +\sup_t \left\{ (1+t)^\frac52  \|\Pp f \|^2 \right\}. \notag
    \end{align}
    Substituting the estimate in Lemma \ref{lemma Pp f} into the above inequality and absorbing the resulting microscopic term, we obtain
    \begin{align}
        &\sup_t\left\{ (1+t)^\frac52  \|w_{\ell_0}\Pp f\|^2 \right\}
        \lesssim  \|w_{{\ell_0}+\frac 32} f_0\|^2          
        +\sup_t\left\{(1+t)^\frac52\|\nabla_x f\|^2\right\}.\notag
    \end{align}
    Finally, applying the decay estimates in Lemmas \ref{lemma:nablaf} and \ref{lemma:k=0} yields the desired estimate \eqref{optimal-decay}.
    This completes the proof of Theorem \ref{thm}.

\section{Appendix}
    For completeness, we present the detailed proof of Lemma \ref{lemma:EE} in the appendix. The lemma establishes the energy estimates used in the global a priori analysis.
    
    \noindent{\bf Proof of Lemma \ref{lemma:EE}.} 
        The proof is divided into three steps.

        \noindent\textit{Step 1. }Energy estimates without weight.
        Applying $\partial^\alpha_x$ with $|\alpha|\leq N$ to \eqref{f-perturbation} and taking the inner product with $\partial^\alpha_xf$ over $\mathbb{R}^3_v\times\mathbb{R}^3_x$, one has
        \begin{align}\label{s1}
            \frac{1}{2}\frac{\rmd}{\rmd t}\|\partial^\alpha_x f\|^2
            + (L\partial^\alpha_xf, \partial^\alpha_xf) 
            =
            ( \partial^\alpha_x \Gamma(f,f), \partial^\alpha_x f), 
        \end{align}
        where we have used the fact that
        \begin{align}
            (v\cdot\nabla_x\partial^\alpha_xf ,\partial^\alpha_x f)=0. \notag
        \end{align}
        For the second term on the left-hand side of \eqref{s1}, it follows from Lemma \ref{lemma L estimate} that
        \begin{align}
            (L\partial^\alpha_xf, \partial^\alpha_xf) \geq \lambda  \|\partial^\alpha_x \Pp f\|^2_\nu. \notag
        \end{align} 
        For the term on the right-hand side of \eqref{s1}, it follows from Lemma \ref{lemma:Gamma-R3} that, for $\alpha_1+\alpha_2=\alpha$, 
        \begin{align}\label{z1}
            |\langle \Gamma^0(\partial^{\alpha_1}_x f, \partial^{\alpha_2}_xf) ,\partial^\alpha_x f \rangle|
			\lesssim &  \sum_{m\leq2} |w_{\ell-m} \nabla_v^m \partial^{\alpha_1}_x  f |_{L^2_v} |\partial_x^{\alpha_2}f|_{L^2_\nu} |\partial_x^{\alpha}f|_{L^2_\nu}, 
		\end{align}
		or
        \begin{align}\label{z2}
            |\langle \Gamma^0(\partial^{\alpha_1}_x f, \partial^{\alpha_2}_xf) ,\partial^\alpha_x f \rangle|
			\lesssim &  \sum_{m\leq2} |w_{\ell-m} \nabla_v^m \partial^{\alpha_2}_x  f |_{L^2_v} |\partial_x^{\alpha_1}f|_{L^2_\nu} |\partial_x^{\alpha}f|_{L^2_\nu}. 
		\end{align}
        If $|\alpha_1|=0$, it follows from \eqref{z1} and the $L^\infty-L^2-L^2$ estimate that 
		\begin{align}
            |( \Gamma^0(\partial^{\alpha_1}_x f,\partial^{\alpha_2}_xf) ,\partial^\alpha_x f )|
			\lesssim& \sum_{m\leq2} \|w_{\ell-m}\nabla_v^m f \|_{L^\infty_x L^2_v}  \|\partial_x^{\alpha}f\|_{L^2_xL^2_\nu}^2 \notag\\
			\lesssim& \sum_{m,n\leq2} \|w_{\ell-m} \nabla_x^n \nabla_v^m f \| \left(  \|\partial_x^{\alpha}\P f \|_\nu^2 +  \|\partial_x^{\alpha} \Pp f \|_\nu^2 \right) \notag\\
			  \lesssim& \sqrt{\mathcal{E}_{N,\ell}(t)} \mathcal{D}_{N,\ell}(t). \notag
		\end{align}
		If $|\alpha_1|=|\alpha|$, the same estimate follows from \eqref{z2} and the $L^\infty-L^2-L^2$ estimate. 
        It remains to consider $1 \leq |\alpha_1| \leq |\alpha|-1$. In this case, $|\alpha_2|+1\leq|\alpha|$.
		If $|\alpha_1|+3 \leq N$, then \eqref{z1} together with the $L^3-L^6-L^2$ estimate gives
		  \begin{align}
            &|( \Gamma^0(\partial^{\alpha_1}_x f,\partial^{\alpha_2}_xf) ,\partial^\alpha_x f )| \notag\\
			\lesssim& \sum_{m\leq 2} \|w_{\ell-m} \nabla_v^m \partial^{\alpha_1}_x  f \|_{L^3_x L^2_v}   \|\partial_x^{\alpha_2}f \|_{L^6_xL^2_\nu}  \|\partial_x^{\alpha}f \|_{L^2_xL^2_\nu}  \notag\\
            \lesssim& \sum_{m\leq 2,n\leq1}  \|w_{\ell-m}\nabla_x^n\nabla_v^m\partial_x^{\alpha_1} f \|  \left(  \|\partial_x^{\alpha_2+e_i} \P f \|_\nu +  \|\partial_x^{\alpha_2+e_i} \Pp f \|_\nu \right) \left(  \|\partial_x^{\alpha} \P f \|_\nu +  \|\partial_x^{\alpha} \Pp f \|_\nu\right) \notag\\
            \lesssim& \sqrt{\mathcal{E}_{N,\ell}(t)} \mathcal{D}_{N,\ell}(t) .\notag
		\end{align}
		If $|\alpha_1|+3 > N$, then $|\alpha_1|\geq N-2$. Since $|\alpha|\leq N$, we have $ |\alpha_2|\leq |\alpha|-N+2\leq2$. Thus, for $N\geq5$, $|\alpha_2|+3\leq5\leq N$. By using \eqref{z2} and the $L^3-L^6-L^2$ estimate, we obtain the same bound.
        For $N=4$, the only remaining case is $|\alpha|=4$ and $|\alpha_1|=|\alpha_2|=2$. 
        Using \eqref{z1} and the $L^2-L^\infty-L^2$ estimate, we obtain
        \begin{align}
            &|( \Gamma^0(\partial_x^{\alpha_1}f,\partial_x^{\alpha_2}f),\partial_x^\alpha f )| \notag\\
            &\lesssim \sum_{m\le 2} \|w_{\ell-m}\nabla_v^m\partial_x^{\alpha_1}f\|_{L^2_xL^2_v}
            \|\partial_x^{\alpha_2}f\|_{L^\infty_xL^2_\nu}
            \|\partial_x^\alpha f\|_{L^2_xL^2_\nu} \notag\\
            &\lesssim \sum_{m\le 2} \|w_{\ell-m}\nabla_v^m\partial_x^{\alpha_1}f\|
            \sum_{n\le 2}\left( \|\nabla_x^n\partial_x^{\alpha_2}\P f\|_\nu + \|\nabla_x^n\partial_x^{\alpha_2}\Pp f\|_\nu\right)
            \left( \|\partial_x^\alpha \P f\|_\nu + \|\partial_x^\alpha \Pp f\|_\nu \right) \notag\\
            &\lesssim \sqrt{\mathcal{E}_{N,\ell}(t)}\,\mathcal{D}_{N,\ell}(t). \notag
        \end{align}
        Plugging the estimates above back into \eqref{s1} and summing over $|\alpha| \leq N$, one has
        \begin{align}\label{s1.1}
            \frac{\rmd}{\rmd t}\sum_{|\alpha|\leq N}\|\partial^\alpha_x f\|^2
            +\lambda \sum_{|\alpha|\leq N}\|\partial^\alpha_x \Pp f\|^2_\nu
            \lesssim 
            \sqrt{\mathcal{E}_{N,\ell}(t)} \mathcal{D}_{N,\ell}(t).
        \end{align}
        Following the strategy in \cite{Duan-2011-Nonlinearity} yields that there exists an interactive functional $|\mathcal{E}_{N,int}(t)|$ satisfying
        \begin{align}
			|\mathcal{E}_{N,int}(t)| &\lesssim \sum_{|\alpha|\leq N} \| \partial^\alpha_x f \|^2, \notag
		\end{align}
        such that for any $t\geq 0$
		\begin{align}
			\frac{\rmd}{\rmd t} \mathcal{E}_{N,int}(t) 
			+ \lambda \sum_{1\leq|\alpha|\leq N} \left|\partial^\alpha_x(a,b,c) \right|^2_{L^2_x} 
			\lesssim 
			\sum_{|\alpha|\leq N}\|\partial^\alpha_x \Pp f \|^2_\nu
			+\sum_{|\alpha|\leq N-1}  \| \langle v \rangle^{-\frac{\gamma}{2}} \partial^\alpha_x\Gamma(f,f)  \|^2. \label{low macro}
		\end{align}
        For the second term on the right-hand side of \eqref{low macro}, it follows from Lemma \ref{lemma:Gamma-L2} that for $\alpha_1+\alpha_2=\alpha$,
        \begin{align}
            |w_\ell\Gamma^0(\partial^{\alpha_1}_xf,\partial^{\alpha_2}_x f)|^2_{L^2_v}
            \lesssim \sum_{m\leq2 }
             | w_{\ell-m} \nabla_v^m \partial^{\alpha_1}_x f |_{L^2_\nu}^2 
             | w_{\ell}\partial^{\alpha_2}_xf |^2_{L^2_\nu},\label{wg1}
        \end{align}
        or
        \begin{align}
            |w_\ell\Gamma^0(\partial^{\alpha_1}_xf,\partial^{\alpha_2}_xf) |^2_{L^2_v}
            \lesssim 
            \sum_{m\leq2 }
             | w_{\ell-m} \nabla_v^m \partial^{\alpha_2}_x f |^2_{L^2_\nu} 
             | w_{\ell}\partial^{\alpha_1}_xf |^2_{L^2_\nu}.\label{wg2}
        \end{align}
        If $|\alpha_1|=0$, it follows from \eqref{wg1} and the $L^\infty-L^\infty-L^2-L^2$ estimate that
        \begin{align}
             \|w_\ell\Gamma^0(\partial^{\alpha_1}_xf,\partial^{\alpha_2}_xf) \|^2
            \lesssim& \sum_{m\leq2 }
             \| w_{\ell-m} \nabla_v^m f \|_{L^\infty_xL^2_v}^2 
             \| w_{\ell}\partial^{\alpha}_xf \|^2_{L^2_xL^2_\nu} \notag\\
            \lesssim& \sum_{m,n\leq2 }
             \| w_{\ell-m} \nabla_x^n \nabla_v^m f \|^2 
            \left(  \| w_{\ell}\partial^{\alpha}_x \P f \|^2_{\nu} +  \| w_{\ell}\partial^{\alpha}_x \Pp f\|^2_{\nu}  \right)\notag\\
            \lesssim& \mathcal{E}_{N,\ell}(t) \mathcal{D}_{N,\ell}(t).\notag
        \end{align}
        If $|\alpha_2|=0$, the same estimate follows from \eqref{wg2} and $L^\infty-L^\infty-L^2-L^2$ estimate. If $1\leq|\alpha_1|\leq|\alpha|-1$ and $|\alpha_1|+3\leq N$, then $|\alpha_2|+1\leq |\alpha|$ and it follows from \eqref{wg1} and $L^3-L^3-L^6-L^6$ estimate that
        \begin{align}
            \|w_\ell\Gamma^0(\partial^{\alpha_1}_xf,\partial^{\alpha_2}_xf) \|^2
            \lesssim& \sum_{m\leq2 }
            \| w_{\ell-m} \nabla_v^m \partial^{\alpha_1}_x f \|_{L^3_xL^2_v}^2
            \| w_{\ell}\partial^{\alpha_2}_xf\|_{L^6_xL^2_\nu}^2
            \notag\\
            \lesssim& \sum_{m\leq2 , n\leq1}
             \| w_{\ell-m}\nabla_x^n \nabla_v^m \partial^{\alpha_1}_x f \|^2 
            \left(  \| w_{\ell}\partial^{\alpha_2+e_i}_x \P f \|^2_{\nu} +  \| w_{\ell}\partial^{\alpha_2+e_i}_x \Pp f\|^2_{\nu}\right) \notag\\
            \lesssim& \mathcal{E}_{N,\ell}(t) \mathcal{D}_{N,\ell}(t).\notag
        \end{align}
        If $1\leq|\alpha_1|\leq|\alpha|-1$ and $|\alpha_1|+3> N$, then $|\alpha_2|+3\leq N$ for $N\geq 5$ and the same estimate follows from \eqref{wg2} and $L^3-L^3-L^6-L^6$ estimate.
        For $N=4$, in the critical case $|\alpha|=4$ and $|\alpha_1|=|\alpha_2|=2$,
        it follows from \eqref{wg1} and $L^2-L^2-L^\infty-L^\infty$ estimate that
        \begin{align}
            \|w_\ell\Gamma^0(\partial^{\alpha_1}_xf,\partial^{\alpha_2}_xf) \|^2
            \lesssim& \sum_{m\leq2 }
            \| w_{\ell-m} \nabla_v^m \partial^{\alpha_1}_x f \|_{L^2_xL^2_v}^2
            \| w_{\ell}\partial^{\alpha_2}_xf\|_{L^\infty_xL^2_\nu}^2
            \notag\\
            \lesssim& \sum_{m\leq2}
             \| w_{\ell-m} \nabla_v^m \partial^{\alpha_1}_x f \|^2 
            \sum_{n\leq2}\left(  \| w_{\ell} \nabla_x^n \partial^{\alpha_2}_x \P f \|^2_{\nu} +  \| w_{\ell} \nabla_x^n \partial^{\alpha_2}_x \Pp f\|^2_{\nu}\right) \notag\\
            \lesssim& \mathcal{E}_{N,\ell}(t) \mathcal{D}_{N,\ell}(t).\notag
        \end{align}
        Plugging the estimates above into \eqref{low macro} and combining with \eqref{s1.1}, one has 
        \begin{align}
            &\frac{\rmd}{\rmd t} \left\{ \sum_{|\alpha|\leq N} \|\partial^\alpha_x f\|^2 + \kappa\mathcal{E}_{N,int}(t) \right\}
            +\lambda \left( \sum_{|\alpha|\leq N} \|\partial^\alpha_x \Pp f \|^2_\nu 
            + \sum_{1 \leq|\alpha|\leq N} |\partial^\alpha_x (a,b,c) |^2_{L^2_x} \right)
            \notag\\
            \lesssim&
            \left( \sqrt{\mathcal{E}_{N,\ell}(t)}+\mathcal{E}_{N,\ell}(t) \right) \mathcal{D}_{N,\ell}(t).\label{ss1.1}
        \end{align}
        \textit{Step 2.} Energy estimates with weight function $w_\ell$.\\
        \noindent\textit{Step 2.1.} Energy estimates without any derivative.
        Applying $\Pp$ to \eqref{f-perturbation} yields the microscopic equation
        \begin{align}\label{s2.1}
            \partial_t\mathbf{P}^{\perp}f + v\cdot\nabla_x\mathbf{P}^\perp f + L\mathbf{P}^{\perp}f = -v\cdot\nabla_x\mathbf{P}f + {\bf P}(v\cdot\nabla_xf) +\Gamma(f,f).
        \end{align}
        Multiplying \eqref{s2.1} by $w_{\ell}^2\Pp f$ and integrating over $\mathbb{R}^3_v\times \mathbb{R}^3_x$, one has
        \begin{align}
            &\frac{1}{2}\frac{\rmd}{\rmd t}\|w_\ell\Pp  f\|^2
            + ( L\Pp f,w_\ell^2\Pp f ) \notag\\
            \lesssim & 
            - ( v\cdot\nabla_x\P f, w_\ell^2\Pp f )
            + ( \P (v\cdot\nabla_x f),w_\ell^2\Pp f )
            + ( \Gamma(f,f), w_\ell^2\Pp f ),  \label{s2.1-1}
        \end{align}
        where we have used the fact that
        \begin{align}
             (v\cdot\nabla_x\Pp f ,w_\ell^2\Pp f )=0.\notag
        \end{align}
        For the second term on the left-hand side of \eqref{s2.1-1}, it follows from Lemma \ref{lemma L estimate} that
        \begin{align}
             ( L\Pp f,w_\ell^2\Pp f ) \geq 
            \lambda \| w_\ell\Pp f \|^2_\nu
            -C\| \chi_{C}\Pp f \|^2.\notag
        \end{align}
        For the first term on the right-hand side of \eqref{s2.1-1}, one has
        \begin{align}
             ( v\cdot\nabla_x\P f, w_\ell^2\Pp f )
            \lesssim& \int_{\mathbb{R}^3_v\times\mathbb{R}^3_x} |v||\nabla_x (a,b,c)| \mu^\delta w_\ell^2 \Pp f\rmd v\rmd x 
            \lesssim \|\nabla_xf\|^2_\nu+\eta\|w_\ell\Pp f \|^2_\nu.\notag
        \end{align}
        For the second term on the right-hand side of \eqref{s2.1-1}, recall the definition of ${\bf P}$ and denote 
		\begin{align}
			\phi_i=\mu^\delta, v_i\mu^\delta, (|v|^2-3)\mu^\delta, \quad i=1,2,3, \notag
		\end{align} 
		then we have 
		\begin{align}
			{\bf P}(v \cdot \nabla_x f)=\int_{\mathbb{R}_v^3} v \cdot \nabla_x f \phi_i \mathrm{d}v \mu^\delta, \notag
		\end{align}
        which implies
		\begin{align}
			 |{\bf P}(v \cdot \nabla_x f) |
			\lesssim \int_{\mathbb{R}_v^3} |v|   |\nabla_x f | |\phi_i| \mathrm{d}v\mu^\delta 
			\lesssim |\nabla_x f|_{L^2_\nu}\mu^\delta. \notag
		\end{align}
		Then one has
        \begin{align}
             ( \P (v \cdot\nabla_xf), w_\ell^2\Pp f )
            \lesssim \int_{\mathbb{R}^3_v\times\mathbb{R}^3_x} |\nabla_x f|_{L^2_\nu} \mu^\delta w_\ell^2 \Pp f\rmd x
            \lesssim \|\nabla_xf\|^2_\nu+\eta\|w_\ell\Pp f \|^2_\nu.\notag
        \end{align}
        For the third term on the right-hand side of \eqref{s2.1-1}, it follows from Lemma \ref{lemma:Gamma-R3} and the $L^\infty-L^2-L^2$ estimate that
        \begin{align}
             ( \Gamma(f,f), w_\ell^2\Pp f )
            \lesssim&  \sum_{m\leq2}\|w_{\ell-m} \nabla_v^m f\|_{L^\infty_xL^2_v}  \|w_\ell f\|_{L^2_xL^2_\nu}\|w_\ell\Pp f\|_{L^2_xL^2_\nu} \notag\\
            \lesssim& \sum_{m,n\leq2}\|w_{\ell-m} \nabla_x^n \nabla_v^m f\| \left( \| w_\ell \P f\|_\nu + \| w_\ell \Pp f\|_\nu \right)\|w_\ell \Pp f\|_\nu \notag\\
            \lesssim& \sqrt{\mathcal{E}_{N,\ell}(t)} \mathcal{D}_{N,\ell}(t).\notag
        \end{align}
        Plugging the estimates above into \eqref{s2.1-1} and choosing a sufficiently small constant $\eta>0$ yields that
        \begin{align}
            \frac{\rmd}{\rmd t}\|w_\ell\Pp  f\|^2
            +\lambda \| w_\ell\Pp f  \|^2_\nu
            \lesssim
            \| \Pp f  \|^2_\nu
            +\|\nabla_xf\|^2_\nu
            +\sqrt{\mathcal{E}_{N,\ell}(t)} \mathcal{D}_{N,\ell}(t).\label{ss2.1}
        \end{align}

        \noindent\textit{Step 2.2.} Weighted estimates on the spatial derivatives.
		Applying $\partial^\alpha_x$ with $1\leq|\alpha|\leq N$ to \eqref{f-perturbation} and taking the inner product with $w_\ell^2\partial^\alpha_xf$ over $\mathbb{R}^3_v\times\mathbb{R}^3_x$, one has
        \begin{align}\label{s2.2}
			\frac{1}{2}\frac{\rmd}{\rmd t}\|w_\ell\partial^\alpha_x f\|^2
            + (L\partial^\alpha_xf,w_\ell^2\partial^\alpha_xf  )
            \lesssim 
             ( \partial^\alpha_x \Gamma(f,f), w_\ell^2\partial^\alpha_x f  ), 
		\end{align}
        where we have used the fact that
        \begin{align}
             (v\cdot\nabla_x\partial^\alpha_xf ,w_\ell^2\partial^\alpha_x f )=0.\notag
        \end{align}
        For the second term on the left-hand side of \eqref{s2.2}, it follows from Lemma \ref{lemma L estimate} that
        \begin{align}
             ( L\partial^\alpha_x f,w_\ell^2\partial^\alpha_x f ) \geq 
            \lambda \| w_\ell\partial^\alpha_x f \|^2_\nu
            -C\| \chi_{C}\partial^\alpha_x f \|^2.\notag
        \end{align}
		For the term on the right-hand side of \eqref{s2.2}, it follows from Lemma \ref{lemma:Gamma-R3} that, for $\alpha_1+\alpha_2=\alpha$, 
		\begin{align}\label{zzz1}
             |\langle \Gamma^0(\partial^{\alpha_1}_x f, \partial^{\alpha_2}_xf), w_\ell^2\partial^\alpha_x f  \rangle|
			\lesssim &  \sum_{m\leq2}   |w_{\ell-m} \nabla_v^m \partial^{\alpha_1}_x f|_{L^2_v}  |w_\ell\partial_x^{\alpha_2}f |_{L^2_\nu}  |w_\ell\partial_x^{\alpha}f |_{L^2_\nu}, 
		\end{align}
		or
		\begin{align}\label{zzz2}
             |\langle \Gamma^0(\partial^{\alpha_1}_x f, \partial^{\alpha_2}_xf), w_\ell^2\partial^\alpha_x f  \rangle|
			\lesssim 
            \sum_{m\leq2}   |w_{\ell-m} \nabla_v^m \partial^{\alpha_2}_x f |_{L^2_v}
             |w_\ell\partial_x^{\alpha_1}f |_{L^2_\nu}
             |w_\ell\partial_x^{\alpha}f |_{L^2_\nu}.
		\end{align}
        If $|\alpha_1|=0$, it follows from \eqref{zzz1} and $L^\infty-L^2-L^2$ estimate that 
		\begin{align}
             |( \Gamma^0(\partial^{\alpha_1}_x f,\partial^{\alpha_2}_xf), w_\ell^2\partial^\alpha_x f  ) |
			\lesssim& \sum_{m\leq2}\|w_{\ell-m} \nabla_v^m f\|_{L^\infty_x L^2_v}  \|w_\ell\partial_x^{\alpha}f\|^2_{L^2_xL^2_\nu} \notag\\
			\lesssim& \sum_{m,n\leq2}\|w_{\ell-m} \nabla_x^n \nabla_v^m f\| \left( \|w_\ell\partial_x^{\alpha} \P f\|^2_\nu +\|w_\ell\partial_x^{\alpha} \Pp f\|^2_\nu \right) \notag\\
			  \lesssim& \sqrt{\mathcal{E}_{N,\ell}(t)} \mathcal{D}_{N,\ell}(t). \notag
		\end{align}
		If $|\alpha_1|=|\alpha|$, the same estimate follows from \eqref{zzz2} and the $L^\infty-L^2-L^2$ estimate. 
		If $1 \leq |\alpha_1| \leq |\alpha|-1$ and $|\alpha_1|+3 \leq N$, then $|\alpha_2|+1\leq|\alpha|$ and it follows from \eqref{zzz1} and the $L^3-L^6-L^2$ estimate that
		  \begin{align}
            &|( \Gamma^0(\partial^{\alpha_1}_x f,\partial^{\alpha_2}_xf) ,w_\ell^2\partial^\alpha_x f)| \notag\\
			\lesssim& \sum_{m\leq2}\|w_{\ell-m}\nabla_v^m\partial^{\alpha_1}_x  f\|_{L^3_x L^2_v}  \|w_\ell\partial_x^{\alpha_2}f\|_{L^6_xL^2_\nu} \|w_\ell\partial_x^{\alpha}f\|_{L^2_xL^2_\nu}  \notag\\
			\lesssim& \sum_{m\leq2,n\leq1}\|w_{\ell-m} \nabla_x^n \nabla_v^m \partial^{\alpha_1}_x  f\|  \left( \|w_\ell\partial_x^{\alpha_2+e_i} \P f\|_\nu +\|w_\ell\partial_x^{\alpha_2+e_i} \Pp f\|_\nu \right)\left( \|w_\ell\partial_x^{\alpha}\P f\|_\nu +\|w_\ell\partial_x^{\alpha} \Pp f\|_\nu \right) \notag\\
            \lesssim& \sqrt{\mathcal{E}_{N,\ell}(t)} \mathcal{D}_{N,\ell}(t). \notag
		\end{align}
		If $1 \leq |\alpha_1| \leq |\alpha|-1$ and $|\alpha_1|+3 > N$, then $|\alpha_2|+3\leq N$ for $N\geq 5$ and the same estimate follows from \eqref{zzz2} and the $L^3-L^6-L^2$ estimate.
        For $N=4$, in the critical case $|\alpha|=4$ and $|\alpha_1|=|\alpha_2|=2$,
        it follows from \eqref{zzz1} and the $L^2-L^\infty-L^2$ estimate that
        \begin{align}
            & ( \Gamma^0(\partial^{\alpha_1}_x f,\partial^{\alpha_2}_xf) ,w_\ell^2\partial^\alpha_x f ) \notag\\
			\lesssim& \sum_{m\leq2}
            \|w_{\ell-m}\nabla_v^m\partial^{\alpha_1}_x  f\|_{L^2_xL^2_v}
            \|w_\ell\partial_x^{\alpha_2}f\|_{L^\infty_xL^2_\nu} \|w_\ell\partial_x^{\alpha}f\|_{L^2_xL^2_\nu}  \notag\\
			\lesssim& \sum_{m\leq2}\|w_{\ell-m} \nabla_v^m \partial^{\alpha_1}_x  f\|  \sum_{n\leq2}\left( \|w_\ell \nabla_x^n \partial_x^{\alpha_2} \P f\|_\nu +\|w_\ell \nabla_x^n  \partial_x^{\alpha_2} \Pp f\|_\nu \right)\left( \|w_\ell\partial_x^{\alpha}\P f\|_\nu +\|w_\ell\partial_x^{\alpha} \Pp f\|_\nu \right) \notag\\
            \lesssim& \sqrt{\mathcal{E}_{N,\ell}(t)} \mathcal{D}_{N,\ell}(t). \notag
		\end{align}
        Plugging the estimates above back into \eqref{s2.2} and summing over $1\leq|\alpha| \leq N$, one has
        \begin{align}\label{ss2.2}
            \frac{\rmd}{\rmd t}\sum_{1\leq|\alpha|\leq N}\|w_\ell\partial^\alpha_x f\|^2
            +\lambda \sum_{1\leq|\alpha|\leq N}\|w_\ell\partial^\alpha_x f\|^2_\nu
            \lesssim 
            \sum_{1\leq|\alpha|\leq N}\|\partial^\alpha_x f\|^2_\nu
            +\sqrt{\mathcal{E}_{N,\ell}(t)} \mathcal{D}_{N,\ell}(t).
        \end{align}
        
        \noindent\textit{Step 2.3.} Weighted estimates on the mixed derivatives. 
        Applying $\partial_x^\alpha \partial_v^\beta$ with $1\leq|\beta|\leq N$ and $|\alpha|+|\beta|\leq N$ to \eqref{s2.1} and taking the inner product with $w_{\ell-|\beta|}^2\partial^\alpha_x\partial^\beta_v\Pp f$ over $\mathbb{R}^3_v\times\mathbb{R}^3_x$, one has
        \begin{align}\label{s2.3}
			&\frac{1}{2}\frac{\rmd}{\rmd t}\|w_{\ell-|\beta|}\partial^\alpha_x\partial^\beta_v \Pp f\|^2
            + (\partial^\alpha_x\partial^\beta_vL\Pp f,w_{\ell-|\beta|}^2\partial^\alpha_x\partial^\beta_v\Pp f  ) \notag\\
            =&
            - ( \partial^\alpha_x\partial^\beta_v(v\cdot\nabla_x \Pp f),w_{\ell-|\beta|}^2\partial^\alpha_x\partial^\beta_v\Pp f )
            - (\partial^\alpha_x\partial^\beta_v(v\cdot\nabla_x\P f), w_{\ell-|\beta|}^2\partial^\alpha_x\partial^\beta_v\Pp f )\notag\\
            &+ (\partial^\alpha_x\partial^\beta_v\P(v\cdot\nabla_xf), w_{\ell-|\beta|}^2\partial^\alpha_x\partial^\beta_v\Pp f )
            + ( \partial^\alpha_x\partial^\beta_v \Gamma(f,f), w_{\ell-|\beta|}^2\partial^\alpha_x\partial^\beta_v\Pp f ). 
		\end{align}
        For the second term on the left-hand side of \eqref{s2.3}, it follows from Lemma \ref{lemma L estimate} that
        \begin{align}
            & (\partial^\alpha_x\partial^\beta_vL\Pp f,w_{\ell-|\beta|}^2\partial^\alpha_x\partial^\beta_v\Pp f  ) \notag\\
            \geq &
            \lambda \| w_{\ell-|\beta|}\partial^\alpha_x \partial^\beta_v \Pp f \|^2_\nu
            -\eta \sum_{|\beta'|\leq |\beta|}\| w_{\ell-|\beta'|}\partial^\alpha_x \partial^{\beta'}_v \Pp f \|^2_\nu
            -C\| \chi_{C}\partial^\alpha_x \Pp f \|^2.\notag
        \end{align}
        For the first term on the right-hand side of \eqref{s2.3}, one has
        \begin{align}
            & ( \partial^\alpha_x\partial^\beta_v(v\cdot\nabla_x \Pp f),w_{\ell-|\beta|}^2\partial^\alpha_x\partial^\beta_v\Pp f )\notag\\
            =&
             ( \partial^{\alpha+e_i}_x\partial^{\beta-e_i}_v\Pp f, w_{\ell-|\beta|}^2\partial^\alpha_x\partial^\beta_v\Pp f  ) 
            +  (v\cdot \nabla_x \partial^\alpha_x\partial^\beta_v\Pp f,w_{\ell-|\beta|}^2\partial^\alpha_x\partial^\beta_v\Pp f )\notag\\
            =&
             (\partial^{\alpha+e_i}_x\partial^{\beta-e_i}_v\Pp f, w_{\ell-|\beta|}^2\partial^\alpha_x\partial^\beta_v\Pp f  )\notag\\
            =&
             (w_{\ell-|\beta-e_i|}\partial^\alpha_x\partial^{\beta-e_i}_v\Pp f\langle v \rangle^\frac{\gamma}{2}, w_{\ell-|\beta|}\partial^\alpha_x\partial^\beta_v\Pp f \langle v \rangle^\frac{\gamma}{2}  )\notag\\
            \lesssim&
            \| w_{\ell-|\beta-e_i|}\partial^\alpha_x\partial^{\beta-e_i}_v\Pp f\|^2_\nu+\eta\| w_{\ell-|\beta|}\partial^\alpha_x\partial^\beta_v\Pp f \|^2_\nu.\notag
        \end{align}
        For the second term on the right-hand side of \eqref{s2.3}, one has
        \begin{align}
             (\partial^\alpha_x\partial^\beta_v(v\cdot\nabla_x\P f), w_{\ell-|\beta|}^2\partial^\alpha_x\partial^\beta_v\Pp f )
            \lesssim&
             ( \partial^{\alpha}_x\partial^\beta_v[\partial^{e_i}_x(a,b,c)\mu^\delta],w_{\ell-|\beta|}^2\partial^\alpha_x\partial^\beta_v\Pp f ) \notag\\
            \lesssim&
             ( \partial^{\alpha+e_i}_x(a,b,c)\mu^{\delta'},w_{\ell-|\beta|}^2\partial^\alpha_x\partial^\beta_v\Pp f ) \notag\\
            \lesssim&
             ( \partial^{\alpha+e_i}_x(a,b,c)\langle v \rangle^\frac{\gamma}{2}\mu^{\delta''},w_{\ell-|\beta|}\partial^\alpha_x\partial^\beta_v\Pp f\langle v \rangle^\frac{\gamma}{2} ) \notag\\
            \lesssim& 
            \| \partial^{\alpha+e_i}_xf \|^2_\nu+\eta\|w_{\ell-|\beta|}\partial^\alpha_x\partial^\beta_v\Pp f\|^2_\nu.\notag
        \end{align}
        For the third term on the right-hand side of \eqref{s2.3}, one has
        \begin{align}
             (\partial^\alpha_x\partial^\beta_v\P(v\cdot\nabla_xf), w_{\ell-|\beta|}^2\partial^\alpha_x\partial^\beta_v\Pp f ) 
            \lesssim&
             (\partial^\alpha_x\partial^\beta_v [\left|\nabla_x f\right|_{L^2_\nu}\mu^\delta], w_{\ell-|\beta|}^2\partial^\alpha_x\partial^\beta_v\Pp f )\notag\\
            \lesssim& 
             (\partial^\alpha_x|\nabla_x f|_{L^2_\nu}\mu^{\delta'},w_{\ell-|\beta|}\partial^\alpha_x\partial^\beta_v\Pp f\langle v \rangle^\frac{\gamma}{2} )\notag\\
            \lesssim& 
            \| \partial^{\alpha+e_i}_xf \|^2_\nu+\eta\|w_{\ell-|\beta|}\partial^\alpha_x\partial^\beta_v\Pp f\|^2_\nu.\notag
        \end{align}
		For the fourth term on the right-hand side of \eqref{s2.3}, it follows from Lemma \ref{lemma:Gamma-R3} that, for $\alpha_1+\alpha_2=\alpha$, $|\beta_1|+|\beta_2|\leq|\beta|$, 
		\begin{align}\label{zz1}
            &|\langle \Gamma^0(\partial^{\alpha_1}_x\partial^{\beta_1}_v f, \partial^{\alpha_2}_x\partial^{\beta_2}_vf), w_{\ell-|\beta|}^2\partial^\alpha_x\partial^{\beta}_v f \rangle| \notag\\
			\lesssim &  
            \sum_{m\leq2} |w_{\ell-m-|\beta_1|} \nabla_v^m \partial^{\alpha_1}_x\partial^{\beta_1}_vf |_{L^2_v} 
             |w_{\ell-|\beta_2|}\partial_x^{\alpha_2}\partial^{\beta_2}_vf |_{L^2_\nu} 
             |w_{\ell-|\beta|}\partial_x^{\alpha}\partial^{\beta}_vf |_{L^2_\nu}, 
		\end{align}
		or
		\begin{align}\label{zz2}
			&|\langle \Gamma^0(\partial^{\alpha_1}_x\partial^{\beta_1}_v f, \partial^{\alpha_2}_x\partial^{\beta_2}_vf), w_{\ell-|\beta|}^2\partial^\alpha_x\partial^{\beta}_v f \rangle| \notag\\
			\lesssim &  
            \sum_{m\leq2}  |w_{\ell-m-|\beta_2|}\nabla_v^m\partial^{\alpha_2}_x\partial^{\beta_2}_vf |_{L^2_v} 
             |w_{\ell-|\beta_1|}\partial_x^{\alpha_1}\partial^{\beta_1}_vf |_{L^2_\nu}  |w_{\ell-|\beta|}\partial_x^{\alpha}\partial^{\beta}_vf |_{L^2_\nu}.
		\end{align}
        If $|\alpha_1|=0$ and $|\beta_1|=0$, it follows from \eqref{zz1} and the $L^\infty-L^2-L^2$ estimate that 
		\begin{align}
			& |( \Gamma^0(\partial^{\alpha_1}_x\partial^{\beta_1}_v f, \partial^{\alpha_2}_x\partial^{\beta_2}_vf) ,w_{\ell-|\beta|}^2\partial^\alpha_x\partial^{\beta}_v f  )|  \notag\\
			\lesssim& 
            \sum_{m\leq2} \|w_{\ell-m} \nabla_v^m f\|_{L^\infty_xL^2_v} 
            \|w_{\ell-|\beta|}\partial_x^{\alpha}\partial^{\beta}_vf\|^2_{L^2_xL^2_\nu} \notag\\
			\lesssim& 
            \sum_{m,n\leq2} \|w_{\ell-m} \nabla_x^n \nabla_v^m f\|  
            \left( \|w_{\ell-|\beta|}\partial_x^{\alpha}\partial^{\beta}_v \P f\|^2_\nu + \|w_{\ell-|\beta|}\partial_x^{\alpha}\partial^{\beta}_v \Pp f\|^2_\nu\right) \notag\\
			  \lesssim& 
            \sqrt{\mathcal{E}_{N,\ell}(t)} \mathcal{D}_{N,\ell}(t). \notag
		\end{align}
		If $|\alpha_1|=|\alpha|$ and $|\beta_1|=|\beta|$, the same estimate follows from \eqref{zz2} and the $L^\infty-L^2-L^2$ estimate. 
		If $1 \leq |\alpha_1|+|\beta_1| \leq |\alpha|+|\beta|-1$ and $|\alpha_1|+|\beta_1|+3 \leq N$, then $|\alpha_2|+|\beta_2|+1\leq |\alpha|+|\beta|$ and it follows from \eqref{zz1} and the $L^3-L^6-L^2$ estimate that
		  \begin{align}
			& |( \Gamma^0(\partial^{\alpha_1}_x\partial^{\beta_1}_v f, \partial^{\alpha_2}_x\partial^{\beta_2}_vf) ,w_{\ell-|\beta|}^2\partial^\alpha_x\partial^{\beta}_v f  )|  \notag\\
			\lesssim& \sum_{m\leq2}\|w_{\ell-m-|\beta_1|}\nabla_v^m\partial^{\alpha_1}_x\partial^{\beta_1}_vf\|_{L^3_x L^2_v}  \|w_{\ell-|\beta_2|}\partial_x^{\alpha_2}\partial^{\beta_2}_vf\|_{L^6_xL^2_\nu} \|w_{\ell-|\beta|}\partial_x^{\alpha}\partial^{\beta}_vf\|_{L^2_xL^2_\nu}  \notag\\
			\lesssim& \sum_{m\leq2,n\leq1}\|w_{\ell-m-|\beta_1|} \nabla_x^n \nabla_v^m \partial^{\alpha_1}_x\partial^{\beta_1}_vf\|  \left( \|w_{\ell-|\beta_2|}\partial_x^{\alpha_2+e_i}\partial^{\beta_2}_v \P f\|_\nu + \|w_{\ell-|\beta_2|}\partial_x^{\alpha_2+e_i}\partial^{\beta_2}_v \Pp f\|_\nu  \right) \notag\\
            &\times\left( \|w_{\ell-|\beta|}\partial_x^{\alpha}\partial^{\beta}_v \P f\|_\nu + \|w_{\ell-|\beta|}\partial_x^{\alpha}\partial^{\beta}_v \Pp f\|_\nu \right) \notag\\
            \lesssim& \sqrt{\mathcal{E}_{N,\ell}(t)} \mathcal{D}_{N,\ell}(t). \notag
		\end{align}
		If $1 \leq |\alpha_1|+|\beta_1| \leq |\alpha|+|\beta|-1$ and $|\alpha_1|+|\beta_1|+3 > N$, then $|\alpha_2|+|\beta_2|+3\leq N$ for $N\geq 5$ and the same estimate follows from \eqref{zz2} and the $L^3-L^6-L^2$ estimate. 
        For $N=4$, in the critical case $|\alpha|+|\beta|=4$ and $|\alpha_1|+|\beta_1|=|\alpha_2|+|\beta_2|=2$,
        it follows from \eqref{zz1} and the $L^2-L^\infty-L^2$ estimate that
        \begin{align}
			& |( \Gamma^0(\partial^{\alpha_1}_x\partial^{\beta_1}_v f, \partial^{\alpha_2}_x\partial^{\beta_2}_vf) ,w_{\ell-|\beta|}^2\partial^\alpha_x\partial^{\beta}_v f  ) | \notag\\
			\lesssim& \sum_{m\leq2}\|w_{\ell-m-|\beta_1|}\nabla_v^m\partial^{\alpha_1}_x\partial^{\beta_1}_vf\|_{L^2_x L^2_v}  \|w_{\ell-|\beta_2|}\partial_x^{\alpha_2}\partial^{\beta_2}_vf\|_{L^\infty_xL^2_\nu} \|w_{\ell-|\beta|}\partial_x^{\alpha}\partial^{\beta}_vf\|_{L^2_xL^2_\nu}  \notag\\
			\lesssim& \sum_{m\leq2}\|w_{\ell-m-|\beta_1|} \nabla_v^m \partial^{\alpha_1}_x\partial^{\beta_1}_vf\|  \sum_{n\leq2}\left( \|w_{\ell-|\beta_2|} \nabla_x^n \partial_x^{\alpha_2}\partial^{\beta_2}_v \P f\|_\nu + \|w_{\ell-|\beta_2|} \nabla_x^n  \partial_x^{\alpha_2}\partial^{\beta_2}_v \Pp f\|_\nu  \right) \notag\\
            &\times\left( \|w_{\ell-|\beta|}\partial_x^{\alpha}\partial^{\beta}_v \P f\|_\nu + \|w_{\ell-|\beta|}\partial_x^{\alpha}\partial^{\beta}_v \Pp f\|_\nu \right) \notag\\
            \lesssim& \sqrt{\mathcal{E}_{N,\ell}(t)} \mathcal{D}_{N,\ell}(t). \notag
		\end{align}
        Plugging the estimates above back into \eqref{s2.3}, summing over $1\leq|\beta|\leq N$ and $|\alpha|+|\beta| \leq N$ and taking a suitable linear combination of the resulting inequalities with constants satisfying $C_1\gg C_2\gg\cdots\gg C_N>0$, one has
        \begin{align}\label{ss2.3}
            &\frac{\rmd}{\rmd t}
            \sum_{|\beta|=1}^N C_{|\beta|} \sum_{ |\alpha|+|\beta| \leq N}\|w_{\ell-|\beta|}\partial^\alpha_x\partial^\beta_v \Pp f\|^2
            +\lambda\sum_{1\leq|\beta|\leq N,|\alpha|+|\beta| \leq N}\| w_{\ell-|\beta|}\partial^\alpha_x \partial^\beta_v \Pp f \|^2_\nu \notag\\
            \lesssim&
            \sum_{|\alpha| \leq N-1}\| w_{\ell}\partial^\alpha_x\Pp f\|^2_\nu +\sum_{1\leq|\alpha|\leq N}\| \partial^{\alpha}_xf \|^2_\nu
            +\sqrt{\mathcal{E}_{N,\ell}(t)} \mathcal{D}_{N,\ell}(t).
        \end{align}
        \textit{Step 3.} Proof of \eqref{energy estimate}. Taking a suitable linear  combination of \eqref{ss1.1}, \eqref{ss2.1}, \eqref{ss2.2} and \eqref{ss2.3}, one has
        \begin{align}
            \frac{\rmd}{\rmd t}\mathcal{E}_{N,\ell} (t)+\lambda \mathcal{D}_{N,\ell}(t) \leq \left(\sqrt{\mathcal{E}_{N,\ell}(t)}+\mathcal{E}_{N,\ell}(t)\right) \mathcal{D}_{N,\ell}(t),\notag
        \end{align}
        where $\mathcal{E}_{N,\ell}(t)$ is defined by
        \begin{align}
            \mathcal{E}_{N,\ell}(t)
            =&
            C\sum_{|\alpha|\leq N} \|\partial^\alpha_x f\|^2 
            +\kappa\mathcal{E}_{N,int}(t)+\|w_\ell\Pp  f\|^2
            +\sum_{1\leq|\alpha|\leq N}\|w_\ell\partial^\alpha_x f\|^2 \notag\\
            &+\sum_{1\leq|\beta|\leq N, |\alpha|+|\beta| \leq N}\|w_{\ell-|\beta|}\partial^\alpha_x\partial^\beta_v \Pp f\|^2 \notag\\
            \sim& \sum_{|\alpha|\leq N} \|\partial^\alpha_x \P f\|^2
            +\sum_{|\alpha|+|\beta| \leq N}\|w_{\ell-|\beta|}\partial^\alpha_x\partial^\beta_v \Pp f\|^2. \notag
        \end{align}
        By choosing $\delta_0$ in \eqref{priori assumption} sufficiently small, we obtain the desired estimate \eqref{energy estimate}. The proof of Lemma \ref{lemma:EE} is complete.

    \bigskip 
    \noindent {\bf Acknowledgements:} 
    The research of Yuanjie Lei  was supported by the National Natural Science Foundation of China (Grant No. 12171117) and the Fundamental Research Funds for the Central Universities. 
    The research of Tong Yang was supported by the General Research Fund of Hong Kong (Project No. 11302723), a start-up grant from The Hong Kong Polytechnic University (Project No. P0043962), and the Kuok Group Foundation. 
    The research of Huijiang Zhao was supported by National Natural Science Foundation of China (Grants 12221001, 12371225 and 12571242).
 
    \vspace{2mm}
    
    \noindent\textbf{Conflict of interest:} The authors do not have any possible conflicts of interest.

    \vspace{2mm}
    
    \noindent\textbf{Data availability statement:}
    Data sharing is not applicable to this article as no data sets were generated or analyzed during the current study.

    \bibliographystyle{plain}

\end{document}